\documentclass[reqno]{amsart}

\usepackage{amsmath}
\usepackage{amsthm}
\usepackage{amscd}
\usepackage{bbm}
\usepackage{enumerate}
\usepackage{amssymb}
\usepackage{palatino}
\usepackage{amscd}
\usepackage{verbatim}
\usepackage{mathrsfs}
\usepackage[margin=1.5in]{geometry}



\title{Self-conjugate core partitions and modular forms}
\author{Levent Alpoge}\email{alpoge@college.harvard.edu}
\address{Quincy House, Harvard College, Cambridge, MA 02138.}

\begin{document}

\begin{abstract}
A recent paper by Hanusa and Nath states many conjectures in the study of self-conjugate core partitions. We prove all but two of these conjectures asymptotically by number-theoretic means. We also obtain exact formulas for the number of self-conjugate $t$-core partitions for ``small'' $t$ via explicit computations with modular forms. For instance, self-conjugate $9$-core partitions are related to counting points on elliptic curves over $\mathbb{Q}$ with conductor dividing $108$, and self-conjugate $6$-core partitions are related to the representations of integers congruent to $11\bmod{24}$ by $3X^2 + 32Y^2 + 96Z^2$, a form with finitely many (conjecturally five) exceptional integers in this arithmetic progression, by an ineffective result of Duke--Schulze-Pillot.
\end{abstract}

\maketitle

\newtheoremstyle{dotless}{}{}{\itshape}{}{\bfseries}{}{ }{}

\newtheorem{thm}{Theorem}
\newtheorem{lem}[thm]{Lemma}
\newtheorem{remark}[thm]{Remark}
\newtheorem{cor}[thm]{Corollary}
\newtheorem{defn}[thm]{Definition}
\newtheorem{prop}[thm]{Proposition}
\newtheorem{conj}[thm]{Conjecture}
\newtheorem{claim}[thm]{Claim}
\newtheorem{exer}[thm]{Exercise}
\newtheorem{fact}[thm]{Fact}

\theoremstyle{dotless}

\newtheorem{thmnodot}[thm]{Theorem}
\newtheorem{lemnodot}[thm]{Lemma}
\newtheorem{cornodot}[thm]{Corollary}

\newcommand{\image}{\mathop{\text{image}}}
\newcommand{\End}{\mathop{\text{End}}}
\newcommand{\Hom}{\mathop{\text{Hom}}}
\newcommand{\Sum}{\displaystyle\sum\limits}
\newcommand{\Prod}{\displaystyle\prod\limits}
\newcommand{\Tr}{\mathop{\mathrm{Tr}}}
\renewcommand{\Re}{\operatorname{\mathfrak{Re}}}
\renewcommand{\Im}{\operatorname{\mathfrak{Im}}}
\newcommand{\im}{\mathrm{im}\,}
\newcommand{\inner}[1]{\langle #1 \rangle}
\newcommand{\pair}[2]{\langle #1, #2\rangle}
\newcommand{\ppair}[2]{\langle\langle #1, #2\rangle\rangle}
\newcommand{\Pair}[2]{\left[#1, #2\right]}
\newcommand{\Char}{\mathop{\mathrm{char}}}
\newcommand{\rank}{\mathop{\mathrm{rank}}}
\newcommand{\sgn}[1]{\mathop{\mathrm{sgn}}(#1)}
\newcommand{\leg}[2]{\left(\frac{#1}{#2}\right)}
\newcommand{\Sym}{\mathrm{Sym}}
\newcommand{\hmat}[2]{\left(\begin{array}{cc} #1 & #2\\ -\bar{#2} & \bar{#1}\end{array}\right)}
\newcommand{\HMat}[2]{\left(\begin{array}{cc} #1 & #2\\ -\overline{#2} & \overline{#1}\end{array}\right)}
\newcommand{\Sin}[1]{\sin{\left(#1\right)}}
\newcommand{\Cos}[1]{\cos{\left(#1\right)}}
\newcommand{\comm}[2]{\left[#1, #2\right]}
\newcommand{\Isom}{\mathop{\mathrm{Isom}}}
\newcommand{\Map}{\mathop{\mathrm{Map}}}
\newcommand{\Bij}{\mathop{\mathrm{Bij}}}
\newcommand{\Z}{\mathbb{Z}}
\newcommand{\R}{\mathbb{R}}
\newcommand{\Q}{\mathbb{Q}}
\newcommand{\C}{\mathbb{C}}
\newcommand{\Nm}{\mathrm{Nm}}
\newcommand{\RI}[1]{\mathcal{O}_{#1}}
\newcommand{\F}{\mathbb{F}}
\renewcommand{\Pr}{\displaystyle\mathop{\mathrm{Pr}}\limits}
\newcommand{\E}{\mathbb{E}}
\newcommand{\coker}{\mathop{\mathrm{coker}}}
\newcommand{\id}{\mathop{\mathrm{id}}}
\newcommand{\Oplus}{\displaystyle\bigoplus\limits}
\renewcommand{\Cap}{\displaystyle\bigcap\limits}
\renewcommand{\Cup}{\displaystyle\bigcup\limits}
\newcommand{\Bil}{\mathop{\mathrm{Bil}}}
\newcommand{\N}{\mathbb{N}}
\newcommand{\Aut}{\mathop{\mathrm{Aut}}}
\newcommand{\ord}{\mathop{\mathrm{ord}}}
\newcommand{\ch}{\mathop{\mathrm{char}}}
\newcommand{\minpoly}{\mathop{\mathrm{minpoly}}}
\newcommand{\Spec}{\mathop{\mathrm{Spec}}}
\newcommand{\Gal}{\mathop{\mathrm{Gal}}}
\newcommand{\Ad}{\mathop{\mathrm{Ad}}}
\newcommand{\Stab}{\mathop{\mathrm{Stab}}}
\newcommand{\Norm}{\mathop{\mathrm{Norm}}}
\newcommand{\Orb}{\mathop{\mathrm{Orb}}}
\newcommand{\pfrak}{\mathfrak{p}}
\newcommand{\qfrak}{\mathfrak{q}}
\newcommand{\mfrak}{\mathfrak{m}}
\newcommand{\Frac}{\mathop{\mathrm{Frac}}}
\newcommand{\Loc}{\mathop{\mathrm{Loc}}}
\newcommand{\Sat}{\mathop{\mathrm{Sat}}}
\newcommand{\inj}{\hookrightarrow}
\newcommand{\surj}{\twoheadrightarrow}
\newcommand{\bij}{\leftrightarrow}
\newcommand{\Ind}{\mathrm{Ind}}
\newcommand{\Supp}{\mathop{\mathrm{Supp}}}
\newcommand{\Ass}{\mathop{\mathrm{Ass}}}
\newcommand{\Ann}{\mathop{\mathrm{Ann}}}
\newcommand{\Krulldim}{\dim_{\mathrm{Kr}}}
\newcommand{\Avg}{\mathop{\mathrm{Avg}}}
\newcommand{\innerhom}{\underline{\Hom}}
\newcommand{\triv}{\mathop{\mathrm{triv}}}
\newcommand{\Res}{\mathrm{Res}}
\newcommand{\eval}{\mathop{\mathrm{eval}}}
\newcommand{\MC}{\mathop{\mathrm{MC}}}
\newcommand{\Fun}{\mathop{\mathrm{Fun}}}
\newcommand{\InvFun}{\mathop{\mathrm{InvFun}}}
\renewcommand{\ch}{\mathop{\mathrm{ch}}}
\newcommand{\irrep}{\mathop{\mathrm{Irr}}}
\newcommand{\len}{\mathop{\mathrm{len}}}
\newcommand{\SL}{\mathrm{SL}}
\newcommand{\GL}{\mathrm{GL}}
\newcommand{\PSL}{\mathrm{SL}}
\newcommand{\actson}{\curvearrowright}
\renewcommand{\H}{\mathbb{H}}
\newcommand{\mat}[4]{\left(\begin{array}{cc} #1 & #2\\ #3 & #4\end{array}\right)}
\newcommand{\interior}{\mathop{\mathrm{int}}}
\newcommand{\floor}[1]{\lfloor #1\rfloor}
\newcommand{\iso}{\cong}
\newcommand{\eps}{\epsilon}
\newcommand{\disc}{\mathrm{disc}}
\newcommand{\Frob}{\mathrm{Frob}}
\newcommand{\charpoly}{\mathrm{charpoly}}
\newcommand{\afrak}{\mathfrak{a}}
\newcommand{\cfrak}{\mathfrak{c}}
\newcommand{\codim}{\mathrm{codim}}
\newcommand{\ffrak}{\mathfrak{f}}
\newcommand{\Pfrak}{\mathfrak{P}}
\newcommand{\homcont}{\hom_{\mathrm{cont}}}
\newcommand{\vol}{\mathrm{vol}}
\newcommand{\ofrak}{\mathfrak{o}}
\newcommand{\A}{\mathbb{A}}
\newcommand{\I}{\mathbb{I}}
\newcommand{\invlim}{\varprojlim}
\newcommand{\dirlim}{\varinjlim}
\renewcommand{\ch}{\mathrm{char}}
\newcommand{\artin}[2]{\left(\frac{#1}{#2}\right)}
\newcommand{\Qfrak}{\mathfrak{Q}}
\newcommand{\ur}[1]{#1^{\mathrm{ur}}}
\newcommand{\absnm}{\mathcal{N}}
\newcommand{\ab}[1]{#1^{\mathrm{ab}}}
\newcommand{\G}{\mathbb{G}}
\newcommand{\dfrak}{\mathfrak{d}}
\newcommand{\Bfrak}{\mathfrak{B}}
\renewcommand{\sgn}{\mathrm{sgn}}
\newcommand{\disjcup}{\bigsqcup}
\newcommand{\zfrak}{\mathfrak{z}}
\renewcommand{\Tr}{\mathrm{Tr}}
\newcommand{\reg}{\mathrm{reg}}
\newcommand{\subgrp}{\leq}
\newcommand{\normal}{\vartriangleleft}
\newcommand{\Dfrak}{\mathfrak{D}}
\newcommand{\nvert}{\nmid}
\newcommand{\K}{\mathbb{K}}
\newcommand{\pt}{\mathrm{pt}}
\newcommand{\RP}{\mathbb{RP}}
\newcommand{\CP}{\mathbb{CP}}
\newcommand{\rk}{\mathop{\mathrm{rk}}}
\newcommand{\redH}{\tilde{H}}
\renewcommand{\H}{\tilde{H}}
\newcommand{\Cyl}{\mathrm{Cyl}}
\newcommand{\T}{\mathbb{T}}
\newcommand{\Ab}{\mathrm{Ab}}
\newcommand{\Vect}{\mathrm{Vect}}
\newcommand{\Top}{\mathrm{Top}}
\newcommand{\Nat}{\mathrm{Nat}}
\newcommand{\inc}{\mathrm{inc}}
\newcommand{\Tor}{\mathrm{Tor}}
\newcommand{\Ext}{\mathrm{Ext}}
\newcommand{\fungrpd}{\pi_{\leq 1}}
\newcommand{\slot}{\mbox{---}}
\newcommand{\funct}{\mathcal}
\newcommand{\Funct}{\mathcal{F}}
\newcommand{\Gunct}{\mathcal{G}}
\newcommand{\FunCat}{\mathrm{Funct}}
\newcommand{\Rep}{\mathrm{Rep}}
\newcommand{\Specm}{\mathrm{Specm}}
\newcommand{\ev}{\mathrm{ev}}
\newcommand{\frpt}[1]{\{#1\}}
\newcommand{\h}{\mathscr{H}}
\newcommand{\poly}{\mathrm{poly}}
\newcommand{\Partial}[1]{\frac{\partial}{\partial #1}}
\newcommand{\Cont}{\mathrm{Cont}}
\renewcommand{\o}{\ofrak}
\newcommand{\bfrak}{\mathfrak{b}}
\newcommand{\Cl}{\mathrm{Cl}}
\newcommand{\ceil}[1]{\lceil #1\rceil}
\newcommand{\hfrak}{\mathfrak{h}}

\let\uglyphi\phi
\let\phi\varphi

\tableofcontents

\section{Introduction}

Since the time of Young it has been known that partitions index the irreducible representations of the symmetric groups. Young and mathematicians of his time also knew that a partition could be encoded in a convenient way --- via what is now known as a \emph{Young diagram} --- and that flipping this diagram about a natural diagonal amounted to tensoring the corresponding irreducible representation with the sign character. Hence it was deduced that the Young diagrams invariant under this flip corresponded to those irreducible representations that split upon restriction to the alternating subgroup.

Some time later, it was discovered by Frame-Robinson-Thrall \cite{framerobinsonthrall} that the hook lengths of a Young diagram determine the dimension of the corresponding irreducible representation (over $\C$). It followed that the study of partitions with hook lengths indivisible by a given integer $t$ --- so-called \emph{$t$-core partitions} --- was connected to modular representation theory.

In this paper we study self-conjugate $t$-core partitions, asymptotically resolving all but two conjectures posed in the paper of Hanusa and Nath \cite{hanusanath} on counting self-conjugate $t$-core partitions. In all but two cases the implied constants are effective, so in principle this reduces many of these conjectures to a finite amount of computation. The ineffective cases are due to the ineffectivity of a result of Duke--Schulze-Pillot \cite{dukeschulzepillot} on integers represented by forms in a given spinor genus, which arises due to the Landau-Siegel phenomenon.

\section{Preliminaries}

Let $\lambda:=\lambda_1\leq\cdots\leq \lambda_k$ be a partition of $n$. For each box $b$ in its associated Young diagram, one defines its hook length $h_b$ by counting the number of boxes directly to its right or below it, including the box itself. The irreducible representations of the symmetric group on $n$ letters, $S_n$, are in explicit bijection with the partitions of $n$. The \emph{hook-length formula} states that the irreducible representation corresponding to $\lambda$ has dimension \begin{align}\dim{\rho_\lambda} = \frac{n!}{\prod h_b},\end{align} the product taken over all the boxes in the Young diagram corresponding to $\lambda$.

The representations of $S_n$ can be defined over $\Z$ (i.e., can be realized as maps $S_n\to \GL_d(\Z)$), and so one may speak of reduction modulo a prime $p$. From modular representation theory one then obtains the criterion that the reduced representation is again irreducible if and only if the general inequality $v_p(n!)\geq v_p(\dim{\rho_\lambda})$ is an equality, where $v_p$ is the $p$-adic valuation. That is, the reduction of $\rho_\lambda$ modulo $p$ is irreducible if and only if none of the $h_b$ are divisible by $p$. This motivates the following more general definition.

\begin{defn}
A partition $\lambda = \lambda_1\leq\cdots \leq \lambda_k$ is called \emph{$t$-core} if none of its hook lengths is divisible by $t$.
\end{defn}

The study of $t$-core partitions goes back at least to Littlewood, who was the first to obtain the generating function for the number of $t$-cores of $n$. Recently Granville and Ono \cite{granvilleono} have resolved precisely which $n$ admit a $t$-core partition, and there has also been activity on a related conjecture of Stanton \cite{stanton} on the monotonicity of $t$-cores in $t$, as well as various identities arising even in Seiberg-Witten theory involving core partitions.

Most relevant to this work is the paper of Hanusa and Nath \cite{hanusanath}, which concerns self-conjugate $t$-cores, or partitions that are both $t$-core and whose Young diagram is symmetric about the natural diagonal (equivalently, those whose corresponding representation $\rho_\lambda$ does \emph{not} remain irreducible upon restriction to the alternating subgroup $A_n\subseteq S_n$). Hanusa and Nath state various conjectures about self-conjugate core partitions, many in direct analogy to conjectures in the study of more general core partitions. In this paper we prove all but two of these conjectures asymptotically.

\section{Main results}

Let $sc_t(n)$ denote the number of self-conjugate $t$-core partitions of $n$, and $sc(n)$ denote the number of self-conjugate core partitions of $n$. By $A\ll_\theta B$ we will mean $|A|\leq C|B|$ for some positive constant $C$ possibly depending on $\theta$. By $A\asymp B$ we will mean $A\gg B$ and $A\ll B$, and by $A\sim B$ we will mean $A = B(1+o(1))$. By $(a,b)$ we will mean the greatest common divisor of $a$ and $b$. For us $\N:=\Z_{\geq 0}$. The greatest integer at most $x$ will be denoted $\floor{x}$. Finally, we will also write $e(z):=e^{2\pi iz}$.

We begin with results on monotonicity. A conjecture of Stanton \cite{stanton} on monotonicity in $t$ of $c_t(n)$, the number of $t$-core partitions of $n$, has a natural analogue for self-conjugate partitions which we can prove asymptotically --- in analogy with Anderson's result \cite{anderson}. (Note that Theorem \ref{monotonicity at four} proves the corresponding result for $t=4$, except with ineffective implied constant.)
\begin{thmnodot}[Cf.\ Conjectures 1.1, 1.2 of Hanusa-Nath \cite{hanusanath}.]\label{self-conjugate monotonicity}
Let $t\geq 9$ or $t=6,8$. Then:
\begin{align}sc_{t+2}(n) > sc_t(n)\end{align} for $n\gg_t 1$, where the implied constant is effectively computable.
\end{thmnodot}

In fact we will prove a slightly more precise result for ``large'' $t$.
\begin{thm}\label{more precise circle method formula} Let $t\geq 10$.
\begin{itemize}
\item If $t$ is even,
\begin{align}
sc_t(n) &= \frac{(2\pi)^{\frac{t}{4}}}{(2t)^{\frac{t}{4}}\Gamma\left(\frac{t}{4}\right)}\left(n + \frac{t^2-1}{24}\right)^{\frac{t}{4}-1}\cdot \left(\sum_{(k,t)=1} k^{-\frac{t}{4}}\sum_{h\in \left(\Z/k\Z\right)^\times} e\left(-\frac{hn}{k}\right)\tilde{\omega}_{h,k}\right) \nonumber\\&\quad\quad+ O_t(n^{\frac{t}{8}})\nonumber\\&\asymp_t \left(n + \frac{t^2-1}{24}\right)^{\frac{t}{4}-1},
\end{align} where $\tilde{\omega}_{h,k}$ is a $24k$-th root of unity, defined precisely in the proof. (See \eqref{definitionofomega}. In fact, the sum over $h$ is a Gauss sum.)
\item If $t$ is odd,
\begin{align}
sc_t(n) &= \frac{(2\pi)^{\frac{t-1}{4}}}{(2t)^{\frac{t-1}{4}}\Gamma\left(\frac{t-1}{4}\right)}\left(n + \frac{t^2-1}{24}\right)^{\frac{t-1}{4}-1}\cdot \left(\sum_{(k,t)=1, k\not\equiv 2\bmod{4}} (2,k)^{\frac{t-1}{4}} k^{-\frac{t-1}{4}}\sum_{h\in \left(\Z/k\Z\right)^\times} e\left(-\frac{hn}{k}\right)\tilde{\omega}_{h,k}\right) \nonumber\\&\quad\quad+ O_t(n^{\frac{t-1}{8}})\nonumber\\&\asymp_t \left(n + \frac{t^2-1}{24}\right)^{\frac{t-1}{4}-1},
\end{align} where $\tilde{\omega}_{h,k}$ is a $24k$-th root of unity, defined precisely in the proof. (Again, the sum over $h$ is a Gauss sum.)
\end{itemize}
\end{thm}

The corresponding result for $c_t(n)$ was proved by Anderson \cite{anderson} using the circle method. Our method is the same for $t\geq 10$, except that we need to be much more explicit in order to bound the leading constants (that is, those in front of the $\left(n + \frac{t^2-1}{24}\right)^\star$ terms --- these are often called the \emph{singular series}) away from $0$. For smaller $t$ we proceed by explicit computation and knowledge of the growth of Fourier coefficients of modular forms. The case $t=4$ will be isolated (see Theorem \ref{monotonicity at four}) due to the ineffectivity of the implied constant.

In the course of the proof of Theorem \ref{self-conjugate monotonicity} we will also establish the following formula.
\begin{thm}\label{exact formula for eight}
For all $n\geq 0$,
\begin{align}
sc_8(n) = \frac{1}{2}\#|\{(x,y,z,w)\in \N^4 : 8n + 21 = x^2 + y^2 + 2z^2 + 2w^2\}|.
\end{align}
As a result, \begin{align} n\ll sc_8(n)\ll n\log\log{n}.\end{align}
\end{thm}

By combining the work of Anderson and Theorem \ref{self-conjugate monotonicity}, we also obtain the following result.
\begin{thmnodot}[Cf.\ Conjecture 4.1 of Hanusa-Nath \cite{hanusanath}.]\label{defect zero blocks and stuff}
Let $11\leq p < q$ be primes. Then the number of defect-zero $p$-blocks of $A_n$ is less than the number of defect-zero $q$-blocks of $A_n$ once $n\gg_{p,q} 1$, where the implied constant is effectively computable in terms of $p$ and $q$.
\end{thmnodot}

Next we move to conjectures about small $t$-cores.

\begin{thmnodot}[Cf.\ Conjecture 3.5 of Hanusa-Nath \cite{hanusanath}.]\label{six cores}
For all $n\geq 0$,
\begin{align}sc_6(n) = \frac{1}{4}\#|\{(x,y,z)\in \N^3 : 24n+35 = 3x^2 + 32y^2 + 96z^2\}|.\end{align}
As a result, $sc_6(n) > 0$ for $n\gg 1$, where the implied constant is \emph{ineffective}.
\end{thmnodot}

Ineffectivity in this paper is due to the Landau-Siegel phenomenon, whereby Siegel's bound $h(-D)\gg_\eps D^{\frac{1}{2}-\eps}$ has ineffective implied constant. (The Goldfeld-Gross-Zagier effective lower bound $h(-D)\gg \log{D}$ is too weak for our purposes.)

\begin{thmnodot}[Cf.\ Conjecture 3.6 of Hanusa-Nath \cite{hanusanath}.]\label{monotonicity at four}
For all $n\geq 0$,
\begin{align}sc_4(n) = \frac{1}{2}\#|\{(x,y)\in \N^2 : 8n+5 = x^2 + y^2\}|.\end{align}
As a result, $sc_6(n) > sc_4(n)$ for $n\gg 1$, where the implied constant is \emph{ineffective}.
\end{thmnodot}

Theorems \ref{six cores} and \ref{monotonicity at four} follow from a computation of the genus and spinor genus of the quadratic form $3X^2 + 32Y^2 + 96Z^2$ (using Magma) and results of Duke--Schulze-Pillot \cite{dukeschulzepillot} (which rely on the subconvexity bound of Iwaniec \cite{iwaniec} for squarefree coefficients of cusp forms of half-integral weight).

Monotonicity is violated for $t=7$, however.
\begin{thmnodot}[Cf.\ Conjecture 3.15 of Hanusa-Nath \cite{hanusanath}.]\label{seven versus nine}
There are infinitely many integers $n$ for which $sc_9(n) < sc_7(n)$.
\end{thmnodot}

The proof is essentially one line: the integers for which $sc_7(n) = 0$ are known (those for which $n+2 = 4^k\cdot (8m+1)$), and, similarly, those for which $sc_9(n) = 0$ are known (those for which $3n+10 = 4^k$). These two sets are infinite and only have $n=2$ in common. The result follows. By computations with Sage, Magma, and Mathematica, we in fact obtain the following more precise formulas for $sc_7(n)$ and $sc_9(n)$.

\begin{thm}\label{exact formula for seven}
For all $n\geq 0$,
\begin{align}
sc_7(n) &= \frac{1}{14}\left(\#|\{(x,y,z)\in \Z^3 : n + 2 = x^2 + y^2 + 2z^2 - yz\}| \right.\nonumber\\&\left.\quad\quad\quad\quad- 2\cdot \#|\{(x,y,z)\in \Z^3 : n+2 = x^2 + 4y^2 + 8z^2 - 4yz\}| \right.\nonumber\\&\left.\quad\quad\quad\quad\quad\quad+ \#|\{(x,y,z)\in \Z^3 : n+2 = 2x^2 + 2y^2 + 3z^2 + 2yz + 2xz + 2xy\}|\right).
\end{align}
\end{thm}

\begin{thm}\label{exact formula for nine}\ \\
\begin{itemize}
\item For $n$ odd,
\begin{align}sc_9(n) = \frac{\sigma(3n+10) + a_{3n+10}(36a) - a_{3n+10}(54a) - a_{3n+10}(108a)}{27}.\end{align}
\item For $n\equiv 0\pmod{4}$,
\begin{align}sc_9(n) = \frac{\sigma(3n+10) + a_{3n+10}(36a) - 3a_{3n+10}(54a) - a_{3n+10}(108a)}{27}.\end{align}
\item For $n\equiv 2\pmod{4}$, writing $3n+10 = 2^e\cdot m$ with $m$ odd,
\begin{align}sc_9(n) = \frac{\sigma(m) + a_{3n+10}(36a) - 3a_{3n+10}(54a) - a_{3n+10}(108a)}{27}.\end{align}
\end{itemize}
Here the $a_n(E)$ are the coefficients appearing in the Dirichlet series for the $L$-function of the elliptic curve $E$. The curve $36a$ is $y^2 = x^3 + 1$, the curve $108a$ is $y^2 = x^3 + 4$, and the curve $54a$ is $y^2 + xy = x^3 - x^2 + 12x + 8$.
\end{thm}

Theorem \ref{exact formula for seven} explains the prevalence of integers congruent to $82\bmod{128}$ appearing in the numerics of Hanusa-Nath \cite{hanusanath}: $3\cdot 82 + 10 = 256$. In fact, looking more closely, the integers $n$ for which $sc_9(n) < sc_7(n)$ that they found all satisfy $3n + 10 = 2^e\cdot m$ with $m$ small and $e$ large.

In any case, by the Hasse bound, we see that, for $n\not\equiv 2\pmod{4}$, \begin{align}sc_9(n) = \frac{\sigma(3n+10)}{27} + O_\eps(n^{\frac{1}{2}+\eps}).\end{align} From this estimate and an elementary construction we see that no inequality of the form $sc_9(n)\gg sc_9(\floor{n/4})$ could possibly hold.

\begin{thmnodot}[Cf.\ Conjecture 4.5 of Hanusa-Nath \cite{hanusanath}.]\label{counterexample}
Let $X > 11$. Let $N_X':=1225\cdot \prod_{7 < p < X} p$, the product taken over the primes between $7$ and $X$. Let $N_X:=N_X'$ (respectively, $2N_X'$) if $N_X'\equiv 1\pmod{3}$ (respectively, $N_X'\equiv 2\pmod{3}$). Let $3n_X+10:=N_X$. Then, for $k = 0, 1, 3, 4$, \begin{align}\frac{sc_9(n_X)}{sc_9(4n_X+k)}\gg \log\log{n_X}.\end{align} Hence, as $X\to\infty$, all parts of Conjecture 4.5 of Hanusa-Nath \cite{hanusanath} are eventually violated.
\end{thmnodot}

Finally, we prove that the proportion of self-conjugate $t$-cores to self-conjugate partitions tends to $1$ if $t$ grows linearly with $n$, in analogy with a result of Craven \cite{craven} to the same effect for $t$-cores proper.

\begin{thmnodot}[Cf.\ Conjecture 4.2 of Hanusa-Nath \cite{hanusanath}.]\label{the one about the proportion}
Let $0 < \alpha < 1$. Then \begin{align}\frac{sc_{\floor{\alpha n}}(n)}{sc(n)}\to 1\end{align} as $n\to\infty$.
\end{thmnodot}

\section{Proofs}

All of the arguments begin from the determination of the generating function for $sc_t(n)$, obtained by Olsson \cite{olsson} and Garvan, Kim, and Stanton \cite{garvankimstanton}. Write \begin{align}F_t(z):=\sum_{n\geq 0} sc_t(n)q^n\end{align} with $q:=e(z)$. Write \begin{align}\eta(z):=q^{\frac{1}{24}}\prod_{n\geq 1} (1 - q^n)\end{align} for the Dedekind eta function.

\begin{thm}\ \\
\begin{itemize}
\item For $t$ even,
\begin{align}q^{\frac{t^2-1}{24}}F_t(z) = \frac{\eta(2z)^2\eta(2tz)^{\frac{t}{2}}}{\eta(z)\eta(4z)}.\end{align}
\item For $t$ odd,
\begin{align}q^{\frac{t^2-1}{24}}F_t(z) = \frac{\eta(2z)^2\eta(2tz)^{\frac{t-5}{2}}\eta(tz)\eta(4tz)}{\eta(z)\eta(4z)}.\end{align}
\end{itemize}
\end{thm}

Hence we see the generating functions are essentially eta products, of weights $\frac{t}{4}$ and $\frac{t-1}{4}$ in the cases of $t$ even and $t$ odd, respectively. These are holomorphic at all cusps, as the following general theorem about eta products (see \cite{kohlerbook}) shows. That the products are holomorphic inside the upper half-plane is immediate from the infinite product representations.

\begin{thm}
Let $f(z):=\prod \eta(mz)^{a_m}$, where $a_m\in \Z$. Then $f$ is holomorphic if and only if for every $c\in \Z$ one has \begin{align}\sum_m \frac{(c,m)^2}{m} a_m\geq 0.\end{align}
\end{thm}

In this case this amounts to the inequalities \begin{align}(c,2)^2 + \frac{(c,2t)^2}{4} - 1 - \frac{(c,4)^2}{4}\geq 0\end{align} for $t$ even, and \begin{align}(c,2)^2 + (c,2t)^2\cdot \frac{t-5}{4t} + \frac{(c,t)^2}{t} + \frac{(c,4t)^2}{4t} - 1 - \frac{(c,4)^2}{4}\geq 0\end{align} for $t$ odd, which both hold by inspection (the expressions are smallest when $c$ and $t$ share no odd prime factor, then split into cases based on $c$ modulo $4$).

We also need explicit formulas for the multiplier systems of the eta and theta functions (with $\theta(z):=\sum_{n\in \Z} q^{n^2}$), which can be found in Knopp \cite{knopp} (except for a missing factor of $2$ in the formula for $v_\theta$), and are originally due to Petersson \cite{petersson}. Before stating the formulas, we set the following notation. For $d$ odd, let \begin{align}\left(\frac{c}{d}\right)_*:=(-1)^{\frac{\sgn(c)-1}{2}\cdot \frac{\sgn(d)-1}{2}}\left(\frac{c}{|d|}\right),\end{align} where $\left(\frac{c}{d}\right)$ is the usual Jacobi symbol, and, for $c$ odd, let \begin{align}\left(\frac{c}{d}\right)^*:=\left(\frac{d}{|c|}\right).\end{align}

\begin{thm}\label{multiplier systems}
Let $\gamma =: \mat{a}{b}{c}{d}\in \SL_2(\Z)$. Then: \begin{align}\eta(\gamma z) = v_\eta(\gamma) (cz+d)^{\frac{1}{2}}\eta(z)\end{align} and, if $\gamma\in \Gamma_0(4)$ (that is to say, $c\equiv 0\bmod{4}$), \begin{align}\theta(\gamma z) = v_\theta(\gamma) (cz+d)^{\frac{1}{2}} \theta(z),\end{align} where we take the principal branch of the square root. More specifically, we have the following formulas for $v_\eta$ and $v_\theta$.
\begin{itemize}
\item The multiplier system of Dedekind's eta function is given by \begin{align}v_\eta(\gamma) = \left(\frac{c}{d}\right)_* e\left(\frac{1}{24}\left((a+d)c - bd(c^2-1) + 3d-3 - 3cd\right)\right)\end{align} if $c$ is even, and \begin{align}v_\eta(\gamma) = \left(\frac{d}{c}\right)^* e\left(\frac{1}{24}\left((a+d)c - bd(c^2-1) - 3c\right)\right)\end{align} if $c$ is odd.
\item The multiplier system of the classical theta function is given by (remember $c\equiv 0\bmod{4}$) \begin{align}v_\theta(\gamma) = \left(\frac{2c}{d}\right)_* e\left(\frac{d-1}{8}\right).\end{align}
\end{itemize}
\end{thm}

With this established, we may begin the arguments.

\subsection{Proof of Theorems \ref{self-conjugate monotonicity}, \ref{more precise circle method formula}, and \ref{exact formula for eight}.}

\subsubsection{Small $t$.}

First, we handle the cases of small $t$ --- the circle method will only tell us something for $t\geq 10$.

We will see that $sc_6(n)$ is proportional to the number of representations of $24n+35$ by the form $3X^2 + 32Y^2 + 96Z^2$. By Siegel's mass formula this is, to leading order, proportional to a class number, which is bounded above by $\ll n^{\frac{1}{2}}\log{n}$. Hence \begin{align}sc_6(n)\ll n^{\frac{1}{2}}\log{n}.\end{align} Since the generating function for $sc_8(n)$ is \begin{align}\sum_{n\geq 0} sc_8(n)q^n = q^{-\frac{21}{8}}\frac{\eta(2z)^2\eta(16z)^4}{\eta(z)\eta(4z)},\end{align} we have that \begin{align}\sum_{n\geq 0} sc_8(n)q^{8n+21} = \left(\frac{\eta(16z)^2}{\eta(8z)}\right)\left(\frac{\eta(64z)^2}{\eta(32z)}\right)\left(\frac{\eta(128z)^2}{\eta(64z)}\right)^2.\end{align} But \begin{align}\frac{\eta(2z)^2}{\eta(z)} = q^{\frac{1}{8}}\sum q^{\triangle_n},\end{align} a shift of the generating function for the triangular numbers. Hence
\begin{align}
\sum_{n\geq 0} sc_8(n)q^{8n+21} &= \left(\sum q^{4n(n+1)+1}\right)\left(\sum q^{16n(n+1)+4}\right)\left(\sum q^{32n(n+1) + 8}\right)^2 \nonumber\\&= \sum_{n\geq 0} \#|\{n = (2a+1)^2 + 4(2b+1)^2 + 8(2c+1)^2 + 8(2d+1)^2, a,b,c,d\geq 0\}|\cdot q^n.
\end{align}

Now $8n+21\equiv 5\pmod{8}$, so that if $8n+21 = x^2 + y^2 + 2z^2 + 8w^2$, without loss of generality we may take $x$ odd and $y$ even. By considering this equality modulo $8$, we see that $4$ does not divide $y$, and hence $2$ divides $z$. Thus the representations of $8n+21$ by the form $X^2 + Y^2 + 2Z^2 + 8W^2$ are equinumerous (modulo switching $X$ and $Y$) with the representations by $X^2 + 4Y^2 + 8Z^2 + 8W^2$. The former is a universal form (as may be easily checked by the Fifteen Theorem \cite{bhargavafifteentheorem}, or looked up in Ramanujan's table of universal diagonal forms \cite{ramanujantableofuniversaldiagonalforms}), and the number of representations of an integer $N$ lies between $\gg N$ and $\ll N\log\log{N}$. Since $sc_8(n)$ is then (up to flipping signs) the number of representations of $8n+21$ by a universal form, we obtain the bounds \begin{align}n\ll sc_8(n)\ll n\log\log{n}.\end{align} (Alternatively, by a result of Shimura \cite{shimura} the theta function of the form is a modular form of weight $2$ and level $8$ with trivial nebentypus, and there are no cusp forms in $M_2(\Gamma_0(8))$.)

The case $t=6$ is proved. We also have the upper bound \begin{align}sc_8(n)\ll_\eps n^{1+\eps}.\end{align} Finally, we will see in Theorem \ref{exact formula for nine} that the same upper bound holds for $sc_9(n)$. Now we will apply the circle method.

\subsubsection{The circle method: even $t$.}

Let \begin{align}P(q):=\frac{q^{\frac{1}{24}}}{\eta(z)},\end{align} the generating function for the partition function $p(n)$. The crux of our calculation is the use of the following transformation formulas. The first is obtained using a transformation formula for the eta function involving a Dedekind sum --- see e.g.\ Apostol \cite{apostol} --- and the second is obtained using a transformation formula for the eta function involving Jacobi symbols --- see e.g.\ Knopp \cite{knopp}.
\begin{thm}
Let $h,k\in \Z$ with $k > 0, (h,k) = 1$. Let $h'\in \Z$ be such that $hh'\equiv -1\pmod{k}$. Let $z\in \C$ be such that $\Re{z} > 0$. Then:
\begin{align}
P\left(e\left(\frac{h}{k} + iz\right)\right) &= e\left(\frac{s(h,k)}{2}\right)\cdot \sqrt{kz}\cdot e^{\frac{\pi}{12k^2z} - \frac{\pi z}{12}}\cdot  P\left(e\left(\frac{h'}{k} + \frac{i}{k^2z}\right)\right),
\end{align}
where \begin{align}s(h,k):=\sum_{r=1}^{k-1} \frac{r}{k}\left(\frac{hr}{k} - \floor{\frac{hr}{k}} - \frac{1}{2}\right)\end{align} is a Dedekind sum. Equivalently,
\begin{align}
P\left(e\left(\frac{h}{k} + iz\right)\right) &= e\left(\frac{1}{24}\left(\frac{h-h'}{k} + iz - \frac{i}{k^2 z}\right)\right)\cdot P\left(e\left(\frac{h'}{k} + \frac{i}{k^2z}\right)\right)\nonumber\\&\quad\quad\quad\quad\cdot \begin{cases} \left(\frac{-h}{k}\right)^* e\left(\frac{1}{24}\left((h'-h)k - \frac{h}{k}(hh'+1)(k^2-1) - 3k\right)\right) & k\text{ odd,}\\ \left(\frac{k}{-h}\right)_* e\left(\frac{1}{24}\left((h'-h)k - \frac{h}{k}(hh'+1)(k^2-1) - 3h - 3 + 3hk\right)\right) & k\text{ even.}\end{cases}
\end{align}
\end{thm}

Write \begin{align}\label{definitionofomega}\omega_{h,k}:=e\left(\frac{s(h,k)}{2}\right),\end{align} a $24k$-th root of unity.

From this calculation we obtain the following transformation formulas for the $F_t$.
\begin{cor} Let $h,k\in \Z$ with $k > 0, (h,k) = 1$.
\begin{itemize}
\item For even $t$, let $h^{(1)},h^{(2)},h^{(3)},h^{(4)}\in \Z$ be such that
\begin{align}
hh^{(1)}&\equiv -1\pmod{k},\\2thh^{(2)}&\equiv -(2t,k)\pmod{k},\\2hh^{(3)}&\equiv -(2,k)\pmod{k},\\\text{and } 4hh^{(4)}&\equiv -(4,k)\pmod{k}.\
\end{align}
Then:
\begin{align}
F_t\left(e\left(\frac{h}{k} + iz\right)\right) &= \frac{\omega_{h,k}\cdot \omega_{\frac{4h}{(4,k)},\frac{k}{(4,k)}}}{\omega_{\frac{2h}{(2,k)},\frac{k}{(2,k)}}^{\ 2}\cdot \omega_{\frac{th}{(2t,k)},\frac{k}{(2t,k)}}^{\frac{t}{2}}}\cdot (kz)^{-\frac{t}{4}}\cdot \sqrt{\frac{(2t,k)^{\frac{t}{2}}\cdot (2,k)^2}{(2t)^{\frac{t}{2}}\cdot (4,k)}}\nonumber\\&\quad\quad\cdot \exp\left(\frac{\pi}{12k^2z}\left(1 + \frac{(4,k)^2}{4} - (2,k)^2 - \frac{(2t,k)^2}{4}\right) + \frac{\pi z}{12}(t^2-1)\right)\nonumber\\&\quad\quad\quad\cdot \frac{P\left(e\left(\frac{h^{(1)}}{k} + \frac{i}{k^2z}\right)\right)\cdot P\left(e\left(\frac{h^{(3)}}{k} + \frac{i(4,k)^2}{4k^2z}\right)\right)}{P\left(e\left(\frac{h^{(4)}}{k} + \frac{i(2,k)^2}{2k^2z}\right)\right)^2\cdot P\left(e\left(\frac{h^{(2)}}{k} + \frac{i(t,k)^2}{2tk^2z}\right)\right)^{\frac{t}{2}}}.
\end{align}

\item For odd $t$, let $h^{(1)},h^{(2)},h^{(3)},h^{(4)},h^{(5)},h^{(6)}\in \Z$ be such that
\begin{align}
hh^{(1)}&\equiv -1\pmod{k},\\4hh^{(2)}&\equiv -(4,k)\pmod{k},\\thh^{(3)}&\equiv -(t,k)\pmod{k},\\4thh^{(4)}&\equiv -(4t,k)\pmod{k},\\2hh^{(5)}&\equiv -(2,k)\pmod{k},\\2thh^{(6)}&\equiv -(2t,k)\pmod{k}.
\end{align}
Then:
\begin{align}
&F_t\left(e\left(\frac{h}{k} + iz\right)\right)\nonumber\\&= \frac{\omega_{h,k}\cdot \omega_{\frac{4h}{(4,k)},\frac{k}{(4,k)}}}{\omega_{\frac{th}{(t,k)},\frac{k}{(t,k)}}\cdot \omega_{\frac{4th}{(4t,k)},\frac{k}{(4t,k)}}\cdot \omega_{\frac{2h}{(2,k)},\frac{k}{(2,k)}}^{\ 2}\cdot \omega_{\frac{2th}{(2t,k)},\frac{k}{(2t,k)}}^{\frac{t-5}{2}}}\nonumber\\&\quad\quad\cdot \sqrt{\frac{k^{-\left(\frac{t-1}{2}\right)}z^{-\left(\frac{t-1}{2}\right)}\cdot (t,k)\cdot (4t,k)\cdot (2,k)^2\cdot (2t,k)^{\frac{t-5}{2}}}{4t^2\cdot (2t)^{\frac{t-5}{2}}\cdot (4,k)}}\nonumber\\&\quad\quad\quad\cdot \exp\left(\frac{\pi}{12k^2z}\left(1 + \frac{(4,k)^2}{4} - \frac{(t,k)^2}{t} - \frac{(4t,k)^2}{4t} - (2,k)^2 - (2t,k)^2\cdot \frac{t-5}{4t}\right) + \frac{\pi z}{12}(t^2-1)\right)\nonumber\\&\quad\quad\quad\quad\cdot \frac{P\left(e\left(\frac{h^{(1)}}{k} + \frac{i}{k^2z}\right)\right)\cdot P\left(e\left(\frac{h^{(2)}}{k} + \frac{i(4,k)^2}{4k^2z}\right)\right)}{P\left(e\left(\frac{h^{(3)}}{k} + \frac{i(t,k)^2}{tk^2z}\right)\right)\cdot P\left(e\left(\frac{h^{(4)}}{k} + \frac{i(4t,k)^2}{4tk^2z}\right)\right)\cdot P\left(e\left(\frac{h^{(5)}}{k} + \frac{i(2,k)^2}{2k^2z}\right)\right)^2\cdot P\left(e\left(\frac{h^{(6)}}{k} + \frac{i(2t,k)^2}{2tk^2z}\right)\right)^{\frac{t-5}{2}}}.
\end{align}
\end{itemize}
\end{cor}

The point of such a formula is to move the argument of $P$ from near the unit circle to near zero (that is, for $|z|$ small), where $P(0) = 1$ gives us total control over the singularities at the roots of unity.

The rest of the calculation follows Anderson rather closely. Let $N\in \Z^+$, and $0 < R < 1$. We take $N\asymp \sqrt{n}$ and $R = e^{-2\pi/n}$. Of course \begin{align}sc_t(n) = \frac{1}{2\pi i}\oint_{\partial\Delta_R} F_t(q) \frac{dq}{q^{n+1}},\end{align} where $\Delta_R\subseteq \C$ is the disk of radius $R$. Now for \begin{align}\frac{h_1}{k_1} < \frac{h}{k} < \frac{h_2}{k_2}\end{align} consecutive Farey fractions of order $N$ (so that $k_1, k, k_2\leq N$), write
\begin{align}
\theta_{h,k}'&:=\frac{1}{k(k+k_1)},\nonumber\\\theta_{h,k}&:=\frac{1}{k(k+k_2)},\nonumber\\\theta_{0,1}', \theta_{0,1}''&:=\frac{1}{N+1}.
\end{align}

The intervals $[\frac{h}{k} - \theta_{h,k}', \frac{h}{k} + \theta_{h,k}'']$ (measure-theoretically) partition $[-\frac{1}{N+1},1 - \frac{1}{N+1}]$, so that \begin{align}sc_t(n) = \sum_{0\leq h < k\leq N, (h,k) = 1} R^{-n}\int_{\frac{h}{k} - \theta_{h,k}'}^{\frac{h}{k} + \theta_{h,k}''} F_t\left(Re(\theta)\right) e(-n\theta) d\theta.\end{align} Writing $R=:e^{-2\pi\eps}$ and $z:=\eps - i\theta$ (we will take $\eps = \frac{1}{n}$), we see that
\begin{align}
sc_t(n) = i\sum_{0\leq h < k\leq N, (h,k) = 1} e\left(-\frac{nh}{k}\right)\int_{\eps + i\theta_{h,k}'}^{\eps - i\theta_{h,k}''} F_t\left(e\left(\frac{h}{k} + iz\right)\right) e^{2\pi nz} dz.
\end{align}

We first do the case of even $t$. In this case, by the transformation formula, we see that
\begin{align}
sc_t(n) &= i\sum_{0\leq h < k\leq N, (h,k) = 1} e\left(-\frac{nh}{k}\right)\cdot \frac{\omega_{h,k}\cdot \omega_{\frac{4h}{(4,k)},\frac{k}{(4,k)}}}{\omega_{\frac{2h}{(2,k)},\frac{k}{(2,k)}}^{\ 2}\cdot \omega_{\frac{2th}{(2t,k)},\frac{k}{(2t,k)}}^{\frac{t}{2}}}\cdot k^{-\frac{t}{4}}\cdot \sqrt{\frac{(2t,k)^{\frac{t}{2}}\cdot (2,k)^2}{(2t)^{\frac{t}{2}}\cdot (4,k)}}\nonumber\\&\quad\quad\quad\quad\cdot \int_{\eps + i\theta_{h,k}'}^{\eps - i\theta_{h,k}''} dz\ z^{-\frac{t}{4}}\cdot \exp\left(\frac{\pi z}{12}\left(24n+t^2-1\right) + \frac{\pi}{12k^2z}\left(1 + \frac{(4,k)^2}{4} - (2,k)^2 - \frac{(2t,k)^2}{4}\right)\right)\nonumber\\&\quad\quad\quad\quad\quad\quad\quad\quad\cdot \frac{P\left(e\left(\frac{h^{(1)}}{k} + \frac{i}{k^2z}\right)\right)\cdot P\left(e\left(\frac{h^{(3)}}{k} + \frac{i(4,k)^2}{4k^2z}\right)\right)}{P\left(e\left(\frac{h^{(4)}}{k} + \frac{i(2,k)^2}{2k^2z}\right)\right)^2\cdot P\left(e\left(\frac{h^{(2)}}{k} + \frac{i(t,k)^2}{2tk^2z}\right)\right)^{\frac{t}{2}}}\nonumber\\&=: M + E_1 + E_2,
\end{align}
where
\begin{align}
M&:=(2t)^{-\frac{t}{4}}i\sum_{0\leq h < k\leq N, (ht,k)=1} e\left(-\frac{nh}{k}\right)\frac{\omega_{h,k}\cdot \omega_{4h,k}}{\omega_{2h,k}^{\ 2}\cdot \omega_{2th,k}^{\frac{t}{2}}}\cdot k^{-\frac{t}{4}}\cdot \int_{\eps + i\theta_{h,k}'}^{\eps - i\theta_{h,k}''} dz\ z^{-\frac{t}{4}} e^{\frac{\pi z}{12}(24n + t^2 - 1)},
\end{align}
\begin{align}
E_1&:=(2t)^{-\frac{t}{4}}i\sum_{0\leq h < k\leq N, (ht,k) = 1} e\left(-\frac{nh}{k}\right)\frac{\omega_{h,k}\cdot \omega_{4h,k}}{\omega_{2h,k}^{\ 2}\cdot \omega_{2th,k}^{\frac{t}{2}}}\cdot k^{-\frac{t}{4}}\nonumber\\&\quad\quad\quad\quad\cdot \int_{\eps + i\theta_{h,k}'}^{\eps - i\theta_{h,k}''} dz\ z^{-\frac{t}{4}}e^{\frac{\pi z}{12}(24n+t^2-1)}\left(\frac{P\left(e\left(\frac{h^{(1)}}{k} + \frac{i}{k^2z}\right)\right)\cdot P\left(e\left(\frac{h^{(3)}}{k} + \frac{i}{4k^2z}\right)\right)}{P\left(e\left(\frac{h^{(4)}}{k} + \frac{i}{2k^2z}\right)\right)^2\cdot P\left(e\left(\frac{h^{(2)}}{k} + \frac{i}{2tk^2z}\right)\right)^{\frac{t}{2}}} - 1\right),
\end{align}
and
\begin{align}
E_2&:=i\sum_{0\leq h < k\leq N, (h,k) = 1, (t,k)\neq 1} e\left(-\frac{nh}{k}\right)\frac{\omega_{h,k}\cdot \omega_{\frac{4h}{(4,k)},\frac{k}{(4,k)}}}{\omega_{\frac{2h}{(2,k)},\frac{k}{(2,k)}}^{\ 2}\cdot \omega_{\frac{2th}{(2t,k)},\frac{k}{(2t,k)}}^{\frac{t}{2}}}\cdot k^{-\frac{t}{4}}\cdot \sqrt{\frac{(2t,k)^{\frac{t}{2}}\cdot (2,k)^2}{(2t)^{\frac{t}{2}}(4,k)}}\nonumber\\&\quad\quad\quad\quad\cdot \int_{\eps + i\theta_{h,k}'}^{\eps - i\theta_{h,k}''} dz\ z^{-\frac{t}{4}}e^{\frac{\pi z}{12}(24n+t^2-1)}\cdot \exp\left(\frac{\pi}{12k^2z}\left(1 + \frac{(4,k)^2}{4} - (2,k)^2 - \frac{(2t,k)^2}{4}\right)\right)\nonumber\\&\quad\quad\quad\quad\quad\quad\quad\quad\cdot \frac{P\left(e\left(\frac{h^{(1)}}{k} + \frac{i}{k^2z}\right)\right)\cdot P\left(e\left(\frac{h^{(3)}}{k} + \frac{i}{4k^2z}\right)\right)}{P\left(e\left(\frac{h^{(4)}}{k} + \frac{i}{2k^2z}\right)\right)^2\cdot P\left(e\left(\frac{h^{(2)}}{k} + \frac{i}{2tk^2z}\right)\right)^{\frac{t}{2}}}.
\end{align}

As suggested by the naming, $M$ will be the main term, and $E_1$ and $E_2$ will be error terms, at least for $t\geq 10$.

Let us first calculate $M$. Note that, on choosing the principal branch of the logarithm on $\C - \R_-$ (the complex plane without the nonpositive reals),
\begin{align}
\int_{\eps + i\theta_{h,k}'}^{\eps - i\theta_{h,k}''} dz\ z^{-\frac{t}{4}}e^{\frac{\pi z}{12}(24n+t^2-1)} &= -\int_{-\infty + i\theta_{h,k}'\to \eps + i\theta_{h,k}'\to \eps - i\theta_{h,k}''\to -\infty - i\theta_{h,k}''} dz\ z^{-\frac{t}{4}} e^{\frac{\pi z}{12}(24n+t^2-1)} \nonumber\\&\quad\quad+ \left(\int_{-\infty + i\theta_{h,k}'\to \eps + i\theta_{h,k}'} - \int_{-\infty - i\theta_{h,k}''\to \eps - i\theta_{h,k}''}\right) dz\ z^{-\frac{t}{4}} e^{\frac{\pi z}{12}(24n+t^2-1)},
\end{align}
where the first integral is over the described contour, with the caveat that the contour does not intersect the nonpositive reals. This path is often called Hankel's contour, since such an integral calculates the gamma function by Hankel's formula. Namely, this becomes
\begin{align}
\int_{\eps + i\theta_{h,k}'}^{\eps - i\theta_{h,k}''} dz\ z^{-\frac{t}{4}}e^{\frac{\pi z}{12}(24n+t^2-1)} &= -\frac{2\pi i}{\Gamma\left(\frac{t}{4}\right)}\left(\frac{\pi(24n+t^2-1)}{12}\right)^{\frac{t}{4}-1} \nonumber\\&\quad\quad+ \left(\int_{-\infty + i\theta_{h,k}'\to \eps + i\theta_{h,k}'} - \int_{-\infty - i\theta_{h,k}''\to \eps - i\theta_{h,k}''}\right) dz\ z^{-\frac{t}{4}} e^{\frac{\pi z}{12}(24n+t^2-1)}.
\end{align}

Hence it suffices to bound these two integrals and the $E_i$.

The integrals pose no problem. Namely, bounding trivially (i.e., via the triangle inequality),
\begin{align}
\int_{-\infty + i\theta_{h,k}'\to \eps + i\theta_{h,k}'} dz\ z^{-\frac{t}{4}} \exp\left(\frac{\pi z}{12}\left(24n + t^2-1\right)\right)&\ll \theta_{h,k}'^{-\frac{t}{4}}\int_{-\infty}^\eps e^{\frac{\pi x}{12}(24n + t^2-1)} dx\nonumber\\&\propto \frac{e^{\frac{\pi\eps}{12}(24n+t^2-1)}}{(24n+t^2-1)\theta_{h,k}'^{\frac{t}{4}}}.
\end{align} (Here $A\propto B$ means $A = cB$ for some constant $c$. That is, $A$ is proportional to $B$.)

Since \begin{align}\theta_{h,k}'^{-1}\gg kN\end{align} by definition, we have the estimate (using $N\asymp \sqrt{n}$ and $\eps^{-1} = n$)
\begin{align}
\int_{-\infty + i\theta_{h,k}'\to \eps + i\theta_{h,k}'} dz\ z^{-\frac{t}{4}} \exp\left(\frac{\pi z}{12}\left(24n + t^2-1\right)\right)\ll_t n^{\frac{t}{8}-1}.
\end{align}

The same holds for the other integral (by Schwarz reflection or repeated effort).

Thus we see that
\begin{align}
M &= \frac{(2\pi)^{\frac{t}{4}}}{(2t)^{\frac{t}{4}}\Gamma\left(\frac{t}{2}\right)}\left(n + \frac{t^2-1}{24}\right)^{\frac{t}{4}-1}\cdot \sum_{0\leq h < k\leq N, (h,k) = 1, (k,t) = 1} e\left(-\frac{nh}{k}\right)\frac{\omega_{h,k}\cdot \omega_{4h,k}}{\omega_{2h,k}^{\ 2}\cdot \omega_{2th,k}^{\frac{t}{2}}}\cdot k^{-\frac{t}{4}}\nonumber\\&\quad\quad + O_t\left(n^{\frac{t}{8}}\right),
\end{align}
where we have used the trivial (and suboptimal) bound \begin{align}\sum_{0\leq h < k\leq N, (h,k) = 1, (k,t) = 1} e\left(-\frac{nh}{k}\right)\frac{\omega_{h,k}\cdot \omega_{4h,k}}{\omega_{2h,k}^{\ 2}\cdot \omega_{2th,k}^{\frac{t}{2}}}\cdot k^{-\frac{t}{4}}\ll n.\end{align}

Next, since all series converge absolutely in the disk, we have the general estimate \begin{align}P\left(e\left(\alpha + \frac{i\beta}{z}\right)\right)^\gamma - 1\ll_\gamma e^{-2\pi\beta\Re\left(1/z\right)},\end{align} obtained by expanding out the relevant series in $q$.

Thus for example
\begin{align}
\frac{P\left(e\left(\frac{h^{(1)}}{k} + \frac{i}{k^2z}\right)\right)\cdot P\left(e\left(\frac{h^{(3)}}{k} + \frac{i}{4k^2z}\right)\right)}{P\left(e\left(\frac{h^{(4)}}{k} + \frac{i}{2k^2z}\right)\right)^2\cdot P\left(e\left(\frac{h^{(2)}}{k} + \frac{i}{2tk^2z}\right)\right)^{\frac{t}{2}}} - 1\ll_t e^{-\frac{\pi\Re\left(1/z\right)}{tk^2}}.
\end{align}

For $z = \eps + iy, y\in [-\theta_{h,k}'',\theta_{h,k}']$, this is \begin{align}\ll_t e^{-\frac{\pi\eps}{tk^2(\eps^2+y^2)}}.\end{align}

Now we turn to $E_2$. We only need that the above bound is $\ll_t 1$. Namely, again bounding trivially,
\begin{align}
E_2&\ll_t \sum_{0\leq h < k\leq N, (h,k) = 1, (k,t)\neq 1} k^{-\frac{t}{4}}\cdot \sqrt{\frac{(2t,k)^{\frac{t}{2}}\cdot (2,k)^2}{(2t)^{\frac{t}{2}}\cdot (4,k)}}\nonumber\\&\quad\quad\quad\quad\cdot \int_{-\theta_{h,k}'}^{\theta_{h,k}''} dy\ (\eps^2 + y^2)^{-\frac{t}{8}} e^{\frac{\pi\eps}{12}(24n + t^2 - 1)}\cdot \exp\left(\frac{\pi\eps}{12k^2(\eps^2+y^2)}\left(1 + \frac{(4,k)^2}{4} - (2,k)^2 - \frac{(2t,k)^2}{4}\right)\right).
\end{align}

Now, if $(t,k)\neq 1$, then (remember $t$ is even!) \begin{align}1 + \frac{(4,k)^2}{4} - (2,k)^2 - \frac{(2t,k)^2}{4}\leq -\frac{1}{4}.\end{align} Also, \begin{align}\sqrt{\frac{(2t,k)^{\frac{t}{2}}\cdot (2,k)^2}{(2t)^{\frac{t}{2}}\cdot (4,k)}}\ll_t 1.\end{align} Hence
\begin{align}
E_2\ll_t \sum_{0\leq h < k\leq N, (h,k) = 1, (t,k)\neq 1} e^{\frac{\pi\eps}{12}(24n + t^2-1)}\cdot \eps^{-\frac{t}{8}}\cdot \int_{-\theta_{h,k}'}^{\theta_{h,k}''} dy\ \left(\frac{\pi\eps}{12k^2(\eps^2+y^2)}\right)^{\frac{t}{8}}\cdot e^{\frac{\pi\eps}{48k^2(\eps^2+y^2)}}.
\end{align}

Note that \begin{align}\frac{\pi\eps}{12k^2(\eps^2+y^2)}e^{-\frac{\pi\eps}{48k^2(\eps^2+y^2)}}\ll 1,\end{align} since the map $x\mapsto xe^{-x}$ is uniformly bounded on $\R_+$. Also, the length of the integral is \begin{align}\theta_{h,k}' + \theta_{h,k}''\ll (kN)^{-1}.\end{align}

Hence
\begin{align}
E_2\ll_t \sum_{0\leq h < k\leq N, (h,k)=1, (t,k)\neq 1} n^{\frac{t}{8}} (kN)^{-1}\ll_t n^{\frac{t}{8}}.
\end{align}

Finally, we turn to bounding $E_1$.

Again bounding trivially (using our ``general bound'')
\begin{align}
E_1&\ll_t \sum_{0\leq h < k\leq N, (h,k) = 1, (t,k) = 1} k^{-\frac{t}{4}}\cdot \int_{-\theta_{h,k}'}^{\theta_{h,k}''} dy\ (\eps^2 + y^2)^{-\frac{t}{8}}\cdot e^{-\frac{\pi\eps}{tk^2(\eps^2+y^2)}}\nonumber\\&\ll_t \sum_{0\leq h < k\leq N, (h,k) = 1, (t,k) = 1} \eps^{-\frac{t}{8}}\cdot (kN)^{-1}\nonumber\\&\ll_t n^{\frac{t}{8}},
\end{align}
as before.

Observing that the difference between the sum with $k\leq N$ and the sum in the theorem statement is
\begin{align}
\sum_{0\leq h < k, k > N, (h,k) = 1, (t,k) = 1} e\left(-\frac{nh}{k}\right) \frac{\omega_{h,k}\cdot \omega_{4h,k}}{\omega_{2h,k}^{\ 2}\cdot \omega_{2th,k}^{\frac{t}{2}}} k^{-\frac{t}{4}}\ll N^{2-\frac{t}{4}}\asymp n^{1 - \frac{t}{8}},
\end{align}
we obtain the first claimed equality for even $t\geq 10$. To show the asymptotic claim, write
\begin{align}
C_t(n) := \sum_{0\leq h < k, (h,k) = 1, (k,t) = 1} e\left(-\frac{nh}{k}\right)\frac{\omega_{h,k}\cdot \omega_{4h,k}}{\omega_{2h,k}^{\ 2}\cdot \omega_{2th,k}^{\frac{t}{2}}}\cdot k^{-\frac{t}{4}}.
\end{align}

Observe that
\begin{align}
|C_t(n) - 1|&\leq \sum_{0 < h < k, (th,k) = 1} k^{-\frac{t}{4}} \nonumber\\&\leq \sum_{k=3,k\text{ odd}}^\infty k^{1 - \frac{t}{4}} \nonumber\\&= \left(1 - 2^{1 - \frac{t}{4}}\right)\zeta\left(\frac{t}{4} - 1\right) - 1\nonumber\\&\leq 0.69.
\end{align}

Thus the asymptotic claim follows, and so we have the full theorem for even $t\geq 10$. We will have to do quite a bit more work in the odd case for $t=11$, but $t\geq 13$ will follow similarly.

\subsubsection{The circle method: odd $t$.}

Things are more complicated in bounding the corresponding $C_t(n)$ for odd $t$, essentially because our eta products do not vanish when $4\vert k$ and $(k,t) = 1$, so there are more terms in the defining sum. But the circle method argument is entirely the same.

The only input is the fact that \begin{align}1 + \frac{(4,k)^2}{4} - \frac{(t,k)^2}{t} - \frac{(4t,k)^2}{4t} - (2,k)^2 - (2t,k)^2\cdot \frac{t-5}{4t} < -1\end{align} if $(t,k)\neq 1$ or $k\equiv 2\bmod{4}$, and it is zero otherwise.

Thus, following the exact same argument as above, we obtain
\begin{align}
sc_t(n) &= \frac{(2\pi)^{\frac{t-1}{4}}}{(2t)^{\frac{t-1}{4}}\Gamma\left(\frac{t-1}{4}\right)}\left(n + \frac{t^2-1}{4}\right)^{\frac{t-1}{4}-1}\nonumber\\&\quad\quad\cdot \sum_{0\leq h < k, (h,k) = 1, (t,k) = 1, k\not\equiv 2\bmod{4}} \left(\frac{(2,k)}{k}\right)^{\frac{t-1}{4}} e\left(-\frac{nh}{k}\right) \frac{\omega_{h,k}\cdot \omega_{\frac{4h}{(4,k)},\frac{k}{(4,k)}}}{\omega_{th,k}\cdot \omega_{\frac{4th}{(4,k)},\frac{k}{(4,k)}}\cdot \omega_{\frac{2h}{(2,k)},\frac{k}{(2,k)}}^{\ 2}\cdot \omega_{\frac{2th}{(2,k)},\frac{k}{(2,k)}}^{\frac{t-5}{4}}}\nonumber\\&\quad\quad\quad\quad + O_t(n^{\frac{t-1}{8}}).
\end{align}

Write, again,
\begin{align}
C_t(n) := \sum_{0\leq h < k, (h,k) = 1, (t,k) = 1, k\not\equiv 2\bmod{4}} \left(\frac{(2,k)}{k}\right)^{\frac{t-1}{4}} e\left(-\frac{nh}{k}\right) \frac{\omega_{h,k}\cdot \omega_{\frac{4h}{(4,k)},\frac{k}{(4,k)}}}{\omega_{th,k}\cdot \omega_{\frac{4th}{(4,k)},\frac{k}{(4,k)}}\cdot \omega_{\frac{2h}{(2,k)},\frac{k}{(2,k)}}^{\ 2}\cdot \omega_{\frac{2th}{(2,k)},\frac{k}{(2,k)}}^{\frac{t-5}{4}}}.
\end{align}

For $t\geq 13$, since $\phi(k)\leq k/2$ for even $k$, we have that \begin{align}|C_t(n) - 1|\leq \zeta\left(\frac{t-1}{4} - 1\right) - 1\leq \zeta(2) - 1 < 0.65.\end{align} Unfortunately a similar argument does not work for $t=11$. So instead we present in the next subsection a calculation that gives \begin{align}|C_{11}(n) - 1|\leq \frac{15609}{854\pi^2} - 1 < 0.86,\end{align} completing the proof.

\subsubsection{Controlling the singular series $C_{11}(n)$.}

Here $t$ will be odd, and soon we will take $t=11$ explicitly.

We will realize the sums over $h$ as Gauss sums. To do this, we will need Petersson's more explicit transformation formula for the eta function, mentioned above (see Theorem \ref{multiplier systems}).

For odd $k$, let $h^{(6)}\in \Z$ be such that $4thh^{(6)}\equiv -1\pmod{k}$, and write
\begin{align}
h^{(1)}&:=4th^{(6)},\\ h^{(2)}&:=th^{(6)},\\ h^{(3)}&:=2th^{(6)},\\ h^{(4)}&:=2h^{(6)},\\ h^{(5)}&:=4h^{(6)},
\end{align}
so that
\begin{align}
hh^{(1)}&\equiv -1\pmod{k},\\ 4hh^{(2)}&\equiv -1\pmod{k},\\ 2hh^{(3)}&\equiv -1\pmod{k},\\ 2thh^{(4)}&\equiv -1\pmod{k},\\ thh^{(5)}&\equiv -1\pmod{k},\\ 4thh^{(6)}&\equiv -1\pmod{k}.
\end{align}

Then Petersson's formula tells us that (after much cancellation --- implicitly we use that $(t,6)=1$, so that $t^2\equiv 1\bmod{24}$, which of course holds in our case)
\begin{align}
F_t\left(e\left(\frac{h}{k} + iz\right)\right) &= (2itkz)^{-\frac{t-1}{4}}\cdot e\left(\frac{1-t^2}{24}\cdot \frac{h}{k}\right)\cdot  e^{\frac{t^2-1}{24}\cdot 2\pi z}\cdot \left(\frac{-2th}{k}\right)^{\frac{t-5}{2}}\cdot e\left(\frac{t-1}{16}k\right)\nonumber\\&\quad\cdot \frac{P\left(e\left(\frac{h^{(1)}}{k} + \frac{i}{k^2z}\right)\right)\cdot P\left(e\left(\frac{h^{(2)}}{k} + \frac{i}{4k^2z}\right)\right)}{P\left(e\left(\frac{h^{(3)}}{k} + \frac{i}{2k^2z}\right)\right)^2\cdot P\left(e\left(\frac{h^{(4)}}{k} + \frac{i}{2tk^2z}\right)\right)^{\frac{t-5}{2}}\cdot P\left(e\left(\frac{h^{(5)}}{k} + \frac{i}{tk^2z}\right)\right)\cdot P\left(e\left(\frac{h^{(6)}}{k} + \frac{i}{4tk^2z}\right)\right)},
\end{align}
where the term $\left(\frac{-2th}{k}\right)$ is the usual Jacobi symbol.

Hence the sum over \emph{odd} $k$ in $C_t(n)$ is, for $t=11$,
\begin{align}
\sum_{k\geq 1, (22,k) = 1} k^{-5/2} e\left(\frac{5k}{8}\right) \left(\frac{-22}{k}\right)\sum_{h\in \Z/k\Z} e\left(-\frac{(n+5)h}{k}\right)\cdot \left(\frac{h}{k}\right),
\end{align}
a Dirichlet series of Gauss sums.

Similarly, in the case of $4\vert k$ (and $(t,k) = 1$), let $h^{(3)}\in \Z$ be such that $thh^{(3)}\equiv -1\pmod{k}$. Write
\begin{align}
h^{(1)}&:=th^{(3)},\\ h^{(2)}&:=h^{(1)},\\ h^{(4)}&:=h^{(3)},\\ h^{(5)}&:=h^{(1)},\\ h^{(6)}&:= h^{(3)},
\end{align}
so that
\begin{align}
hh^{(1)}&\equiv -1\pmod{k},\\ 4hh^{(2)}&\equiv -4\pmod{k},\\ thh^{(3)}&\equiv -1\pmod{k},\\ 4thh^{(4)}&\equiv -4\pmod{k},\\ 2hh^{(5)}&\equiv -2\pmod{k},\\ 2thh^{(6)}&\equiv -2\pmod{k}.
\end{align}

Then, applying our transformation formulas with these $h^{(i)}$, after a great deal of cancellation we see that, for $t=11$,
\begin{align}
F_t\left(e\left(\frac{h}{k} + iz\right)\right) &= -e\left(-\frac{5h}{k}\right)\cdot e\left(-5iz\right)\cdot \left(\frac{2k}{h}\right)\cdot e\left(\frac{5h}{8}\right)\cdot e\left(\frac{5}{8}\right)\cdot e\left(\frac{hk}{16}\right)\cdot \begin{cases} \left(\frac{11}{k/4}\right) & 8\nmid k\\ \left(\frac{k/4}{11}\right) e\left(\frac{h}{4}\right)\cdot e\left(-\frac{hk}{16}\right) & 8\vert k\end{cases}.
\end{align}

Thus the sum over even $k$ in $C_{11}(n)$ can be written (splitting into a sum over odd $k$ and $e\geq 2$ via replacing $k$ by $2^e k$)
\begin{align}
-e\left(\frac{5}{8}\right)&\left(\sum_{k\geq 1, (22,k)=1} \left(\frac{11}{k}\right)\cdot (2k)^{-5/2}\cdot \sum_{h\in \Z/4k\Z} \left(\frac{8k}{h}\right)\cdot e\left(\frac{k(5+k/2)-2(n+5)}{8k}\cdot h\right) \right.\nonumber\\&\quad\quad\left.+ \sum_{e > 2}\sum_{k\geq 1, (22,k)=1} (-1)^e\cdot \left(\frac{k}{11}\right)\cdot (2^{e-1}k)^{-5/2}\cdot \sum_{h\in \Z/2^ek\Z} \left(\frac{2^{e+1}k}{h}\right)\cdot e\left(\frac{-(n+5)-2^{e-3}k}{2^ek}\cdot h\right)\right).
\end{align}

We can evaluate Gauss sums (or, perhaps more correctly, ``twisted Ramanujan sums'') exactly (see Montgomery-Vaughan \cite{montgomeryvaughan} Theorem 9.12).
\begin{thm}
Let $\chi$ be a Dirichlet character of conductor $d\vert q$, and let $\chi^*$ be the corresponding primitive character inducing $\chi$. Then:
\begin{align}
\sum_{a\in \Z/q\Z} \chi(a)e\left(\frac{an}{q}\right) = \begin{cases} 0 & d\nmid \frac{q}{(n,q)}\\ \overline{\chi^*}\left(\frac{n}{(q,n)}\right)\cdot \chi^*\left(\frac{q}{(n,q)d}\right)\cdot \mu\left(\frac{q}{(n,q)d}\right)\cdot \frac{\phi(q)}{\phi\left(\frac{q}{(n,q)}\right)}\cdot \tau(\chi^*) & \text{otherwise,}\end{cases}
\end{align}
where $\tau(\chi^*)$ is the Gauss sum corresponding to $\chi^*$, of absolute value $\sqrt{d}$ if $\chi^*$ is nonprincipal, and $\mu$ is the usual Mobius function.
\end{thm}

So, for odd $k$, since $\left(\frac{\cdot}{k}\right)$ is primitive modulo the squarefree part of $k$ (which we will denote $k/\square$, where $\square$ is the largest square dividing $k$), we have that
\begin{align}
&\sum_{h\in \Z/k\Z} e\left(-\frac{(n+5)h}{k}\right)\cdot \left(\frac{h}{k}\right) \nonumber\\&\quad\quad= \begin{cases} 0 & (n+5,k)\nmid \square\\ \left(\frac{-(n+5)/(n+5,k)}{k/\square}\right)\cdot \left(\frac{\square/(n+5,k)}{k/\square}\right)\cdot \mu\left(\frac{\square}{(n+5,k)}\right)\cdot \frac{\phi(k)}{\phi\left(\frac{k}{(n+5,k)}\right)}\cdot \tau\left(\left(\frac{\cdot}{k/\square}\right)\right) & \text{otherwise.}\end{cases}
\end{align}

Next we turn to the even Gauss sums. Since $\left(\frac{8k}{\cdot}\right)$ is primitive modulo $8\left(\frac{k}{\square}\right)$ (again $k$ is odd), we see that, writing $\gcd:=(n+5,k)$,
\begin{align}
&\sum_{h\in \Z/4k\Z} \left(\frac{8k}{h}\right)\cdot e\left(\frac{k(5+k/2)-2(n+5)}{8k}\cdot h\right) \nonumber\\&\quad\quad= \frac{1}{2}\sum_{h\in \Z/8k\Z} \left(\frac{8k}{h}\right)\cdot e\left(\frac{k(5+k/2)-2(n+5)}{8k}\cdot h\right) \nonumber\\&\quad\quad= \frac{1}{2}\begin{cases} 0 & \gcd\nmid \square,\\ \left(\frac{8(k/\square)}{(k(5+k/2)-2(n+5))/\square}\right)\cdot \left(\frac{8(k/\square)}{\square/\gcd}\right)\cdot \mu\left(\frac{\square}{\gcd}\right)\cdot \frac{\phi(8k)}{\phi\left(\frac{8k}{\gcd}\right)}\cdot \tau\left(\left(\frac{8(k/\square)}{\cdot}\right)\right) & \text{otherwise,}\end{cases}
\end{align} where the first equality follows from considering $h\mapsto h+4k$ in $\Z/8k\Z$ --- the summand picks up a minus sign from each term, and so does not change.

Similarly, for $e > 2$, since $\left(\frac{2^{e+1}k}{\cdot}\right)$ is primitive modulo $\left(\frac{k}{\square}\right)\cdot 2^r := \left(\frac{k}{\square}\right)\cdot \begin{cases} 1 & e\text{ odd, } k\equiv 1\pmod{4},\\ 4 & e\text{ odd, } k\equiv 3\pmod{4},\\ 8 & e\text{ even,}\end{cases}$ we see that, writing $\gcd:=(n+5 + 2^{e-3}k,2^ek)$,
\begin{align}
&\sum_{h\in \Z/2^ek\Z} \left(\frac{2^{e+1}k}{h}\right)\cdot e\left(\frac{-(n+5)-2^{e-3}k}{2^ek}\cdot h\right) \nonumber\\&\quad\quad= \begin{cases} 0 & \gcd\nmid 2^{e-r}\cdot \square,\\ \left(\frac{2^r(k/\square)}{(-(n+5)-2^{e-3}k)/\gcd}\right)\cdot \left(\frac{2^r(k/\square)}{2^{e-r}\square/\gcd}\right)\cdot \mu\left(\frac{2^{e-r}\square}{\gcd}\right)\cdot \frac{\phi(2^ek)}{\phi\left(\frac{2^ek}{\gcd}\right)}\cdot \tau\left(\left(\frac{2^r (k/\square)}{\cdot}\right)\right) & \text{otherwise.}\end{cases}
\end{align}

As horrible and unweildy as these formulas may look, the essential observation is that their \emph{absolute values} are (almost) multiplicative in $k$ (that is, the absolute value of the term corresponding to $k\ell$ is the product of those corresponding to $k$ and $\ell$ if $(k,\ell) = 1$). Namely, for $k$ odd, writing $k =: \prod_{v_p(k)\text{ odd}} p\cdot \prod_p p^{e_p}$ with each $e_p\in 2\Z$ and $v_p(\cdot)$ the $p$-adic valuation (so that the second term is precisely what we have been calling $\square$) and $\gcd:=(n+5,k)$, we have the following formulas.

First,
\begin{align}
&\left|\sum_{h\in \Z/k\Z} \left(\frac{h}{k}\right)\cdot e\left(-\frac{(n+5)h}{k}\right)\right| \nonumber\\\quad&= \begin{cases} 0 & \gcd\nmid \square,\text{ or }v_p(\gcd)\neq e_p\text{ if $v_p(k)$ odd, or }v_p(\gcd) < e_p-1\text{ if $v_p(k)$ even,}\\ \frac{\phi(k)}{\phi\left(\frac{k}{\gcd}\right)}\sqrt{\frac{k}{\square}} & \text{otherwise.}\end{cases}
\end{align}
This is multiplicative in $k$.

Next, if $e > 2$ and $e$ is \emph{odd} ($\gcd = (n+5,k)$ still), then
\begin{align}
&\left|\sum_{h\in \Z/2^ek\Z} \left(\frac{2^{e+1}k}{h}\right)\cdot e\left(\frac{-(n+5)-2^{e-3}k}{2^ek}\cdot h\right)\right| \nonumber\\&\quad= \begin{cases} 0 & \gcd\nmid \square,\text{ or }v_p(\gcd)\neq e_p\text{ if $v_p(k)$ odd},\\&\quad\text{ or }v_p(\gcd) < e_p-1\text{ if $v_p(k)$ even},\\&\quad\quad\text{ or, if $k\equiv 1\bmod{4}$, }2^{e-1}\nmid (n+5)+2^{e-3}k,\\&\quad\quad\quad\quad\text{ or, if $k\equiv 3\bmod{4}$, }v_2(n+5)\neq e-2\\ 2^{e-1}\frac{\phi(k)}{\phi\left(\frac{k}{\gcd}\right)}\sqrt{\frac{k}{\square}\cdot \begin{cases} 1 & k\equiv 1\bmod{4}\\ 4 & k\equiv 3\bmod{4}\end{cases}} & \text{otherwise.}\end{cases}
\end{align}
This is \emph{not} multiplicative in $k$, but, by weakening conditions on being zero a bit and factoring out the terms depending only on $e$, we can bound it above by something that is. Namely,
\begin{align}
&\left|\sum_{h\in \Z/2^ek\Z} \left(\frac{2^{e+1}k}{h}\right)\cdot e\left(\frac{-(n+5)-2^{e-3}k}{2^ek}\cdot h\right)\right|\nonumber\\&\quad\leq \begin{cases} 0 & v_2(n+5)\not\in \{e-2, e-3\},\text{ or }\gcd\nmid \square,\\&\quad\text{ or }v_p(\gcd)\neq e_p\text{ if $v_p(k)$ odd,}\\&\quad\quad\text{ or }v_p(\gcd) < e_p-1\text{ if $v_p(k)$ even}\\ 2^{e-1}\frac{\phi(k)}{\phi\left(\frac{k}{\gcd}\right)}\sqrt{\frac{k}{\square}\cdot \begin{cases} 1 & v_2(n+5) = e-3\\ 4 & v_2(n+5) = e-2\end{cases}} & \text{otherwise.}\end{cases}
\end{align}

If $e > 2$ is even,
\begin{align}
\left|\sum_{h\in \Z/2^ek\Z} \left(\frac{2^{e+1}k}{h}\right)\cdot e\left(\frac{-(n+5)-2^{e-3}k}{2^ek}\right)\right| = \begin{cases} 0 & 2^{e-2}\nmid n+5,\text{ or } \gcd\nmid\square,\\&\quad\text{ or } v_p(\gcd)\neq e_p\text{ if $v_p(k)$ odd},\\&\quad\quad\text{ or } v_p(\gcd) < e_p - 1\text{ if $v_p(k)$ even}\\ 2^{e - \frac{3}{2}}\frac{\phi(k)}{\phi\left(\frac{k}{\gcd}\right)}\sqrt{\frac{k}{\square}} & \text{otherwise.}\end{cases}
\end{align}
This is multiplicative in $k$ once we factor out the terms depending only on $e$.

Finally, for $e=2$,
\begin{align}
\left|\sum_{h\in \Z/4k\Z} \left(\frac{8k}{h}\right)\cdot e\left(\frac{k(5+k/2) - 2(n+5)}{8k}\right)\right| = \begin{cases} 0 & \gcd\nmid\square,\\&\quad\text{ or } v_p(\gcd)\neq e_p\text{ if $v_p(k)$ odd},\\&\quad\quad\text{ or } v_p(\gcd) < e_p - 1\text{ if $v_p(k)$ even}\\ \sqrt{2}\frac{\phi(k)}{\phi\left(\frac{k}{\gcd}\right)}\sqrt{\frac{k}{\square}} & \text{otherwise.}\end{cases}
\end{align}
Note that the right-hand side is the same result as setting $e=2$ in the $e>2$, $e$ even formula. In particular this is also multiplicative in $k$ once we factor out the terms depending only on $e$.

The formulas may look horrendous, but we are about to apply them for prime powers only (thanks to multiplicativity), where they become rather simple.

For instance, the sum over odd $k$ (so $e=0$) in $|C_t(n) - 1|$ is at most
\begin{align}
&\sum_{k > 1, k\text{ odd, }(k,11) = 1} k^{-5/2}\left|\sum_{h\in \Z/k\Z} \left(\frac{h}{k}\right)\cdot e\left(-\frac{(n+5)h}{k}\right)\right|\nonumber\\&\quad\quad= \sum_{k > 1, k\text{ odd, }(k,11) = 1} k^{-5/2}\begin{cases} 0 & (n+5,k)\nmid \square,\\&\quad\text{ or }v_p(\gcd)\neq e_p\text{ if $v_p(k)$ odd},\\&\quad\quad\text{ or }v_p(\gcd) < e_p-1\text{ if $v_p(k)$ even,}\\ \frac{\phi(k)}{\phi\left(\frac{k}{(n+5,k)}\right)}\sqrt{\frac{k}{\square}} & \text{otherwise}\end{cases}
\nonumber\\&\quad\quad = \left(\prod_{p\neq 2,11, v_p(n+5)\text{ odd}} \sum_{a=0}^{\frac{v_p(n+5)+1}{2}} p^{-5a}\frac{\phi(p^{2a})}{\begin{cases} 1 & a < \frac{v_p(n+5)+1}{2}\\ \phi(p) & \text{otherwise}\end{cases}}\right.\nonumber\\&\quad\quad\quad\quad\cdot \left.\prod_{p\neq 2,11, v_p(n+5)\text{ even}} \left(\sum_{a=0}^{\frac{v_p(n+5)}{2}} p^{-5a}\phi(p^{2a}) + p^{-\frac{5}{2}(v_p(n+5)+1)}\frac{\phi\left(p^{v_p(n+5)+1}\right)\sqrt{p}}{\phi(p)}\right)\right) - 1.
\end{align}
This ends up simplifying to \begin{align}\left(\prod_{p\neq 2,11} \frac{1-p^{-4}}{1-p^{-3}}\left(1 + \frac{p^{-2-3\floor{\frac{v_p(n+5)}{2}}}}{1+p}\right)\right) - 1 =: D(n) - 1.\end{align}

The same holds for the other sums, too. That is, the sum over $e > 0$ is bounded above by
\begin{align}
&\sum_{e=2, e\text{ even}}^{v_p(n+5)+2} (2^{e-1})^{-5/2}\sum_{k\geq 1, (22,k)=1} k^{-5/2}\left(2^{e-\frac{3}{2}}\begin{cases} 0 & \star\\ \frac{\phi(k)}{\phi\left(\frac{k}{(n+5,k)}\right)}\sqrt{\frac{k}{\square}} & \text{otherwise}\end{cases}\right) \nonumber\\&\quad\quad+ \begin{cases} (2^{(v_2(n+5)+2)-1})^{-5/2}\sum_{k\geq 1, (22,k)=1} k^{-5/2} \left(2^{v_2(n+5)+2}\begin{cases} 0 & \star\\ \frac{\phi(k)}{\phi\left(\frac{k}{(n+5,k)}\right)}\sqrt{\frac{k}{\square}} & \text{otherwise}\end{cases}\right) & v_2(n+5)\text{ odd,}\\ (2^{(v_2(n+5)+3)-1})^{-5/2}\sum_{k\geq 1, (22,k)=1} k^{-5/2} \left(2^{(v_2(n+5)+3)-1}\begin{cases} 0 & \star\\ \frac{\phi(k)}{\phi\left(\frac{k}{(n+5,k)}\right)}\sqrt{\frac{k}{\square}} & \text{otherwise}\end{cases}\right) & v_2(n+5)\text{ even,}\end{cases}
\end{align}
where by the condition ``$\star$'' we mean: $(n+5,k)\nmid \square,\text{ or }v_p(\gcd)\neq e_p\text{ if $v_p(k)$ odd, or }v_p(\gcd) < e_p-1\text{ if $v_p(k)$ even.}$

That is, it is bounded above by
\begin{align}
\sum_{e\geq 1}^{\floor{\frac{v_p(n+5)+2}{2}}} 2^{1-3e} D(n) + \begin{cases} 2^{-\frac{3v_2(n+5)+1}{2}} D(n) & v_2(n+5)\text{ odd,}\\ 2^{-\frac{3v_2(n+5)+6}{2}} D(n) & v_2(n+5)\text{ even.}\end{cases}
\end{align}

So, adding up the odd and even contributions (and subtracting $1$ from the odd sum), we get:
\begin{align}
|C_{11}(n) - 1|&\leq D(n)\cdot \left(1 + 2\sum_{e = 1}^{\floor{\frac{v_2(n+5)+2}{2}}} 2^{-3e} + \begin{cases} 2^{-\frac{3v_2(n+5)+1}{2}} & v_2(n+5)\text{ odd}\\ 2^{-\frac{3v_2(n+5)+6}{2}} & v_2(n+5)\text{ even}\end{cases}\right) - 1\\&\leq D(n)\left(\frac{9}{7} + \frac{1}{4}\right) - 1.
\end{align}

Since $D(n)\leq \prod_{p\neq 2,11} (1 + p^{-2}) = \frac{\zeta(2)}{\zeta(4)}\cdot \left(\frac{(1-2^{-4})(1-11^{-4})}{(1-2^{-2})(1-11^{-2})}\right)$, we see that \begin{align}|C_{11}(n) - 1|\leq \frac{15609}{854\pi^2} - 1 = 0.8519\ldots,\end{align} as desired. This completes the proofs of Theorems \ref{self-conjugate monotonicity} and \ref{more precise circle method formula}.

\subsection{Proof of Theorem \ref{defect zero blocks and stuff}.}

As noted in Hanusa-Nath \cite{hanusanath}, the number of defect-zero $p$-blocks of $A_n$ is \begin{align}2sc_p(n) + \frac{1}{2}\left(c_p(n) - sc_p(n)\right) = \frac{1}{2}c_p(n) + \frac{3}{2}sc_p(n).\end{align} By work of Anderson \cite{anderson}, \begin{align}c_p(n)\asymp_p n^{\frac{p-3}{2}},\end{align} and as we saw above \begin{align}sc_p(n)\ll_p n^{\frac{p-5}{4}}.\end{align} The result follows.

\subsection{Proof of Theorem \ref{six cores}.}

The generating function for $6$-cores is
\begin{align}
\sum_{n\geq 0} sc_6(n) q^{24n+35} &= \left(\frac{\eta(48z)^2}{\eta(24z)}\right)\left(\frac{\eta(288z)^3}{\eta(96z)}\right) \nonumber\\&= \left(\sum_{n\geq 0} q^{3(2n+1)^2}\right)\left(\frac{\eta(288z)^3}{\eta(96z)}\right).
\end{align}
The second factor is the generating function \begin{align}\sum_{n\geq 0} c_3(n)q^{96n + 32},\end{align} where $c_3(n)$ denotes the number of $3$-cores of $n$. By an identity of Jacobi (see e.g.\ \cite{onoonthreecores}) the coefficients are known: \begin{align}c_3(n) = \sum_{d\vert 3n+1} \left(\frac{d}{3}\right).\end{align} Note that the right-hand side is a multiplicative function of $3n+1$. On prime powers $p^k$ with $p\equiv 1\pmod{3}$ it takes the value $k+1$, and on prime powers $p^k$ with $p\equiv 2\pmod{3}$ it takes the values $0$ or $1$ according to whether $k$ is odd or even, respectively. By classical algebraic number theory, this is exactly half of the number of representations of $3n+1$ by the form $X^2 + 3Y^2$. That is, \begin{align}c_3(n) = \frac{1}{2}\#|\{3n+1 = x^2 + 3y^2\}|,\end{align} or
\begin{align}
\sum_{n\geq 0} sc_6(n) q^{24n+35} &= \frac{1}{2}\sum_{n\geq 0} \#|\{n = 3(2a+1)^2 + 32b^2 + 96c^2, a\geq 0\}|\cdot q^n.
\end{align}
Note that if $24n+35 = 3x^2 + 32y^2 + 96z^2$, then $x$ must be odd.

Hence we obtain the claimed formula. A computation in Magma shows that the spinor genus of the form $3X^2 + 32Y^2 + 96Z^2$ coincides with its genus. By a theorem of Duke--Schulze-Pillot \cite{dukeschulzepillot} this gives the ineffective claim about positivity of $sc_6(n)$, since the integers $24n+35$ are locally represented by this form.

\subsection{Proof of Theorem \ref{monotonicity at four}.}

We have already seen that (ineffectively) \begin{align}sc_6(n)\gg_\eps n^{\frac{1}{2}-\eps}.\end{align} The generating function for $4$-cores is
\begin{align}
\sum_{n\geq 0} sc_4(n)q^{8n+5} &= \left(\frac{\eta(16z)^2}{\eta(8z)}\right)\left(\frac{\eta(64z)^2}{\eta(32z)}\right) \nonumber\\&= \left(\sum_{n\geq 0} q^{(2n+1)^2}\right)\left(\sum_{n\geq 0} q^{4(2n+1)^2}\right) \nonumber\\&= \sum_{n\geq 0} \#|\{n = (2a+1)^2 + 4(2b+1)^2, a,b\geq 0\}|\cdot q^n.
\end{align}
Of course if $8n+5 = a^2 + b^2$, without loss of generality $a$ is odd and $b\equiv 2\pmod{4}$, so we see that \begin{align}sc_4(n) = \frac{1}{4}\#|\{8n+5 = x^2 + y^2\}|.\end{align} Writing $8n + 5 =: \prod_p p^{e_p}$, we know that the right-hand side is precisely \begin{align}\prod_{p\equiv 1\bmod{4}} (e_p+1),\end{align} or $0$ if there is a $p\equiv 3\bmod{4}$ with $e_p$ odd.

In either case this is bounded above by \begin{align}\prod_{p\equiv 1\bmod{4}} (e_p+1)\leq \tau(n)\ll_\eps n^\eps,\end{align} whence the monotonicity result for $t=4$.

\subsection{Proof of Theorem \ref{seven versus nine}.}

This was proved in the Main Results section (see the statement of Theorem \ref{seven versus nine}): let $n:=\frac{4^k-10}{3}$ for $k>2$. Then $sc_9(n) = 0$, but $n+2 = 4\cdot \left(\frac{4^{k-1}-1}{3}\right)$, which is $4$ times something congruent to $5\bmod{8}$, so $sc_7(n)\neq 0$.

\subsection{Proof of Theorem \ref{exact formula for seven}.}

Let \begin{align}G_t(z):=q^{-2}F_7(z).\end{align} Then, for $\gamma=:\mat{a}{b}{c}{d}\in \Gamma_0(28)$, a calculation with the multplier systems for the eta and theta functions shows that \begin{align}G_t(\gamma z) = \left(\frac{7}{d}\right) \left(\frac{\theta(\gamma z)}{\theta(z)}\right)^3 G_t(z).\end{align} Hence $G_t$ is a modular form of weight $3/2$ of level $28$ and nebentypus character $\chi_7:=\left(\frac{7}{d}\right) = (-1)^{\frac{d-1}{2}}\left(\frac{d}{7}\right)$.

A paper of Lehman \cite{lehmanmathcomputation} lists the ternary quadratic forms of level $28$ and discriminant $7$ --- by a theorem of Shimura \cite{shimura} these forms have associated theta functions of weight $3/2$, level $28$, and nebentypus $\chi_7$ as well. The forms are \begin{align}X^2 + Y^2 + 2Z^2 - XZ,\end{align}\begin{align}X^2 + Y^2 + 7Z^2 \text{ (which is in the same genus as $X^2 + Y^2 + 4Z^2 - 2YZ$)},\end{align}\begin{align}2X^2 + 2Y^2 + 3Z^2 + 2YZ + 2XZ + 2XY,\end{align}\begin{align}\text{and } X^2 + 4Y^2 + 8Z^2 - 4YZ.\end{align}

A computation in Sage shows that the theta functions associated to these quadratic forms form a basis for the four-dimensional space of modular forms of weight $3/2$, level $28$, and nebentypus $\chi_7$. Using Sage to express $G_t$ in terms of this basis gives the claimed formula.

\begin{remark}
In fact a finite computation in Sage does amount to a proof of the equality for all $n$, since we can easily check that $G_t$ and the sum above have $q$-expansions agreeing well past the Sturm bound, which is smaller than $100$ in all cases. Hence, since both sides are modular, they must agree for all $n$ (the point is that the space is finite-dimensional). Note that the same remark applies for the following subsection as well.
\end{remark}

\subsection{Proof of Theorem \ref{exact formula for nine}.}

Let \begin{align}H_t(z):=\sum_{n\geq 0} sc_9(n) q^{3n+10} = \frac{\eta(6z)^2\eta(54z)^2\eta(27z)\eta(108z)}{\eta(3z)\eta(12z)}.\end{align} A calculation with the multiplier system for the eta function shows that $H_t$ is a modular form of weight $2$, level $108$, and trivial nebentypus character.

According to Sage,
\begin{align}
H_t(z) &= -\frac{2}{27}E_{\chi_3,\chi_3}(q^4) + \frac{1}{54}E_{\chi_3,\chi_3}(q) + \frac{2}{81}\tilde{E}_{1,1}(q^3) \nonumber\\&\quad\quad+ \frac{1}{54}\tilde{E}_{1,1}(q^4) - \frac{1}{162}\tilde{E}_{1,1}(q^9) - \frac{2}{81}\tilde{E}_{1,1}(q^{12}) + \frac{1}{162}\tilde{E}_{1,1}(q^{36}) \nonumber\\&\quad\quad+ f(z),
\end{align} where $f(z)$ is a cusp form, $\chi_3 := \left(\frac{\cdot}{3}\right)$ is the Legendre symbol modulo $3$, $\tilde{E}_{1,1}(q^t):=E_{1,1}(q) - tE_{1,1}(q^t)$, and \begin{align}E_{\chi,\psi}(q):=\sum_{n\geq 0}\left(\sum_{d\vert n} \chi(d)\psi(n/d) d\right) q^n.\end{align}

Another computation in Sage gives us that \begin{align}f(z) = \frac{1}{27}f_{36a}(q) - \frac{1}{54}\left(f_{54a}(q) + f_{54a}^{\chi_3}(q)\right) + \frac{1}{27}\left(f_{54a}(q^2) - f_{54a}^{\chi_3}(q^2)\right) - \frac{1}{27}f_{108a}(q),\end{align} where $f_{\cdots}$ is the eigenform associated to the elliptic curve ``$\cdots$'', and $f_{54a}^{\chi_3}$ indicates twisting by $\chi_3$ (so $f_{54a}^{\chi_3} = f_{54b}$).

Extracting coefficients gives the result.

\subsection{Proof of Theorem \ref{counterexample}.}

By Theorem \ref{exact formula for nine}, we see that the quotient \begin{align}\frac{sc_9(n)}{sc_9(4n+k)} = \frac{\sigma(3n+10)}{\sigma(12n+3k+10)} + O_\eps(n^{-\frac{1}{2}+\eps}).\end{align} Now, writing $N_X=:2^a\cdot 5^2\cdot 7^2\cdot N_X'$ (so that $a=0$ or $1$), \begin{align}\sigma(N_X) = 1767\cdot \prod_{7 < p < X} (p+1) = N_X\cdot \exp\left(\sum_{7 < p < X} \frac{1}{p} + O(1)\right)\gg N_X\log{X}.\end{align} Also, by the prime number theorem (recall $3n_X + 10 = N_X$), \begin{align}\log{n_X}\sim \log{N_X}\sim X,\end{align} so that this lower bound is \begin{align}\sigma(N_X)\gg N_X\log\log{n_X}.\end{align} Next, $4N_X + 3k - 30$ is not divisible by any prime $7 < p < X$, since $N_X$ is and $3(10-k)$ is not (recall $0\leq k\leq 4$). Write \begin{align}4N_X + 3k - 30 =: 2^b\cdot 5^c\cdot 7^d\cdot \prod_{p > X} p^{e_p}.\end{align} Note that $b,c,d\leq 1$ since $k = 0,1,3,4$. Thus
\begin{align}
\sigma(4N_X + 3k - 30) &= (4N_X + 3k - 30)\cdot (2 - 2^{-b})\cdot \left(\frac{5-5^{-c}}{4}\right)\cdot \left(\frac{7-7^{-d}}{6}\right)\cdot \prod_{p > X} \left(1 + p^{-1} + \cdots + p^{-e_p}\right) \nonumber\\&= (4N_X + 3k - 30)\cdot \exp\left(\sum_{p > X: e_p\neq 0} \frac{1}{p} + O(1)\right).
\end{align}
Certainly \begin{align}\sum_{p > X: e_p\neq 0} \frac{1}{p}\leq \sum_{X < p\leq Y} \frac{1}{p}\end{align} for $Y$ chosen so that $\pi(Y) - \pi(X) = \#|\{p : e_p\neq 0\}|$. Note that \begin{align}\prod_{X < p\leq Y} p\leq \prod_{p : e_p\neq 0} p\ll N_X,\end{align} whence, by the prime number theorem, \begin{align}Y - X\leq X + o(X).\end{align} Thus the sum \begin{align}\sum_{X < p\leq Y} \frac{1}{p} = \log\log{Y} - \log\log{X} + O(1)\ll \log\log(3X) - \log\log(X) + O(1)\ll 1\end{align} remains bounded.

Hence \begin{align}\sigma(12n_X + 3k + 10) = \sigma(4N_X + 3k - 30)\ll N_X,\end{align} proving the claim.

\subsection{Proof of Theorem \ref{the one about the proportion}.}

The following argument is in exact analogy with that of Craven \cite{craven} for $c_t(n)$.

In Hanusa-Nath \cite{hanusanath}, the following theorem is established.

\begin{thmnodot}[Theorems 3.4 and 3.11 of Hanusa-Nath \cite{hanusanath}.]\label{recursion theorem from hanusa nath}
Let $\hat{p}_t(x)$ be the number of ordered sequences of $t$ partitions $\lambda_1, \ldots, \lambda_t$ of integers $a_1,\ldots,a_t$ such that $\sum_{i=1}^t a_i = x$. (Note that this is polynomial in $t$ for $t\gg_x 1$.)
\begin{itemize}
\item The number of self-conjugate $2t$-cores is given by \begin{align}sc_{2t}(n) = \sum_I (-1)^k \hat{p}_t(i_1)\cdots \hat{p}_t(i_a) sc(n - 4tk),\end{align} the sum taken over $I = (i_1, \ldots, i_a)$ (each $i_j > 0$) such that $k := i_1 + \cdots + i_a\leq \frac{n}{4t}$.
\item The number of self-conjugate $(2t+1)$-cores is given by \begin{align}sc_{2t+1}(n) = \sum_{I,J} (-1)^k \hat{p}_t(i_1)\cdots \hat{p}_t(i_a) sc(j_1)\cdots sc(j_a) sc(n - (2t+1)(2k+\ell)),\end{align} the sum taken over $I = (i_1, \ldots, i_k)$ and $J = (j_1, \ldots, j_k)$ (each $i_a, j_a\geq 0$ and $i_a + j_a > 0$) such that $2k+\ell := 2(i_1 + \cdots + i_a) + (j_1 + \cdots + j_a)\leq \frac{n}{2t+1}$.
\end{itemize}
\end{thmnodot}

For the proof of Theorem \ref{the one about the proportion} we only need the following corollary.
\begin{cor} In the same notation as Theorem \ref{recursion theorem from hanusa nath}, the number of self-conjugate partitions which are \emph{not} $t$-cores satisfies
\begin{align}sc(n) - sc_t(n)\ll \sum_{I,J : I\neq \emptyset} \hat{p}_t(i_1)\cdots \hat{p}_t(i_k) sc(j_1)\cdots sc(j_k) sc(n - t(2k+\ell)).\end{align}
\end{cor}

Thus 
\begin{align}\label{proportion error term}
1 - \frac{sc_{\floor{\alpha n}(n)}}{sc(n)}\ll \sum_{I,J : I\neq \emptyset} \hat{p}_{\floor{\alpha n}}(i_1)\cdots \hat{p}_{\floor{\alpha n}}(i_k) sc(j_1)\cdots sc(j_k) \frac{sc(n - \floor{\alpha n}(2k+\ell))}{sc(n)},
\end{align}
where the sum over $I,J$ is a sum over $O(1)$ terms --- specifically, it is over $I,J$ such that $I\neq \emptyset$ and for which \begin{align}2k + \ell\leq \frac{n}{\floor{\alpha n}}\ll_\alpha 1.\end{align}

Now the important observation is that \begin{align}\log{sc(n)}\sim c\sqrt{n}\end{align} for some constant $c$, which can be proved elementarily from the formula \begin{align}\sum_{n\geq 0} sc(n)q^n = \prod_{n\geq 0} (1 + q^{2n+1}),\end{align} as in Hardy-Ramanujan \cite{hardyramanujanonp(n)}. Hence \begin{align}\log{sc(n)} - \log{sc(n - \floor{\alpha n}(2k+\ell))}\gg \sqrt{n}\end{align} since $2k+\ell > 0$. Thus each term \begin{align}\frac{sc(n - \floor{\alpha n}(2k+\ell))}{sc(n)}\end{align} of \eqref{proportion error term} decays faster than any polynomial, whence the result.

\section{Acknowledgements}

This research was conducted at the University of Minnesota Duluth REU program, supported by NSF/DMS grant 1062709 and NSA grant H98230-11-1-0224. I would like to thank Joe Gallian for his constant encouragement and for the wonderful environment for research at UMD. I would also like to thank Will Jagy and Noam Elkies for greatly helpful discussions related to ternary quadratic forms with more than two exceptional integers represented by their genus.

\bibliography{selfconjugatecorepartitionsandmodularforms}{}

\begin{thebibliography}{10}

\bibitem{anderson}
Jaclyn Anderson.
\newblock An asymptotic formula for the {$t$}-core partition function and a
  conjecture of {Stanton}.
\newblock {\em Journal of Number Theory}, 128:2591--2615, 2008.

\bibitem{apostol}
Tom~M. Apostol.
\newblock {\em Modular functions and {Dirichlet} series in number theory}.
\newblock Springer-Verlag, New York, 1990.

\bibitem{bhargavafifteentheorem}
Manjul Bhargava.
\newblock On the {Conway-Schneeberger} fifteen theorem.
\newblock In {\em Quadratic forms and their applications}, volume 272 of {\em
  Contemporary Mathematics}, pages 27--37, Providence, Rhode Island, 1999.
  American Mathematical Society.

\bibitem{craven}
David~A. Craven.
\newblock The number of {$t$}-cores of size {$n$}.
\newblock See http://people.maths.ox.ac.uk/craven/docs/papers/tcores0608.pdf,
  unpublished (2006).

\bibitem{dukeschulzepillot}
William Duke and Rainer Schulze-Pillot.
\newblock Representations of integers by positive ternary quadratic forms and
  equidistribution of lattice points on ellipsoids.
\newblock {\em Inventiones Mathematicae}, 99:49--57, 1990.

\bibitem{framerobinsonthrall}
J.~S. Frame, G.~de~B. Robinson, and R.~M. Thrall.
\newblock The hook graphs of the symmetric groups.
\newblock {\em Canadian J. Math.}, 6:316--324, 1954.

\bibitem{garvankimstanton}
Frank Garvan, Dongsu Kim, and Dennis Stanton.
\newblock Cranks and {$t$}-cores.
\newblock {\em Inventiones Mathematicae}, 101(1):1--17, 1990.

\bibitem{granvilleono}
Andrew Granville and Ken Ono.
\newblock Defect zero {$p$}-blocks for finite simple groups.
\newblock {\em Transactions of the American Mathematical Society},
  348:331--347, 1996.

\bibitem{onoonthreecores}
Guo-Niu Han and Ken Ono.
\newblock Hook lengths and 3-cores.
\newblock {\em Annals of Combinatorics}, 15:305--312, 2011.

\bibitem{hanusanath}
Christopher~R.H. Hanusa and Rishi Nath.
\newblock The number of self-conjugate core partitions.
\newblock {\em Journal of Number Theory}, 133:751--768, 2013.

\bibitem{hardyramanujanonp(n)}
Godfrey~Harold Hardy and Srinivasa Ramanujan.
\newblock Asymptotic formulae in combinatory analysis.
\newblock {\em Proceedings of the London Mathematical Society}, pages 75--115,
  1918.

\bibitem{iwaniec}
Henryk Iwaniec.
\newblock Fourier coefficients of modular forms of half integral weight.
\newblock {\em Inventiones Mathematicae}, 87:385--401, 1987.

\bibitem{knopp}
Marvin~I. Knopp.
\newblock {\em Modular functions in analytic number theory}.
\newblock AMS Chelsea Publishing, Providence, Rhode Island, 1993.

\bibitem{kohlerbook}
Gunter Kohler.
\newblock {\em Eta products and theta series identities: examples}.
\newblock Springer, Berlin, 2011.

\bibitem{lehmanmathcomputation}
J.~Larry Lehman.
\newblock Levels of positive definite ternary quadratic forms.
\newblock {\em Mathematics of Computation}, 58:399--417, 1992.

\bibitem{montgomeryvaughan}
Hugh~L. Montgomery and Robert~C. Vaughan.
\newblock {\em Multiplicative number theory I: classical theory}.
\newblock Cambridge studies in advanced mathematics, Cambridge, UK, 2007.

\bibitem{olsson}
Jorn~B. Olsson.
\newblock On the {$p$}-blocks of symmetric and alternating groups and their
  covering groups.
\newblock {\em Journal of Algebra}, 128:188--213, 1990.

\bibitem{petersson}
Hans Petersson.
\newblock Uber die arithmetischen eigenschaften eines systems multiplikativer
  modulfunktionen von primzahlstufe.
\newblock {\em Acta Mathematica}, 95:57--110, 1956.

\bibitem{ramanujantableofuniversaldiagonalforms}
Srinivasa Ramanujan.
\newblock On the expression of a number in the form {$ax^2 + by^2 + cz^2 +
  du^2$}.
\newblock {\em Proceedings of the Cambridge Philosophical Society}, 19:11--21,
  1917.

\bibitem{shimura}
Goro Shimura.
\newblock On modular forms of half integral weight.
\newblock {\em Annals of Mathematics}, 97(3):440--481, 1973.

\bibitem{stanton}
Dennis Stanton.
\newblock Open positivity conjectures for integer partitions.
\newblock {\em Trends in Mathematics}, 2:19--25, 1990.

\end{thebibliography}
\bibliographystyle{plain}

\ \\

\end{document}